\documentclass{article}
\usepackage[utf8]{inputenc}
\usepackage[margin=1in, includehead, includefoot]{geometry}
\usepackage{amsthm,amssymb,amsmath}
\usepackage{enumitem}
\usepackage{mathrsfs}
\usepackage{mathtools}
\usepackage{fancyhdr}
\usepackage{indentfirst}
\usepackage{graphicx}
\usepackage{placeins}
\usepackage{wrapfig}
\usepackage{url}

\newtheorem{definition}{Definition}
\newtheorem{lemma}{Lemma}[section]
\newtheorem{theorem}{Theorem}[section]
\newtheorem{corollary}{Corollary}[section]
\newtheorem{observation}{Observation}[section]
\newtheorem{prop}{Proposition}[section]
\newtheorem{conjecture}{Conjecture}

\DeclarePairedDelimiter\ceil{\lceil}{\rceil}

\newcommand{\toroman}[1]{\textit{\expandafter{\romannumeral #1\relax}}}

\newcommand{\cbeginproof}[0]{\par\noindent\textit{Proof.} }
\newcommand{\cendproof}[0]{ \qed\par\vspace{1em}}

\newcommand{\npcompleteproblem}[3]{ \par\vspace{0.5em}\noindent{\textbf{#1}}\newline \textbf{INSTANCE: } #2 \newline \textbf{QUESTION: } #3 \par\vspace{0.5em} }

\newcommand{\centered}[1]{\begin{tabular}{c} #1 \end{tabular}}

\title{Optimal Error-detection system for Identifying Codes}
\author{
    \small Devin C. Jean\\
    \small Computer Science Department \\
    \small Vanderbilt University\\
    \small \texttt{devin.c.jean@vanderbilt.edu}
    \and
    \small Suk J. Seo\\
    \small Computer Science Department\\
    \small Middle Tennessee State University\\
    \small \texttt{Suk.Seo@mtsu.edu}
}
\date{}

\begin{document}
\maketitle
\thispagestyle{empty}

\begin{abstract}
Assume that a graph $G$ models a detection system for a facility with a possible ``intruder," or a multiprocessor network with a possible malfunctioning processor.
We consider the problem of placing detectors at a subset of vertices in $G$ to determine the location of an intruder if there is any.
Many types of detection systems have been defined for different sensor capabilities; in particular, we focus on Identifying Codes, where each detector can determine whether there is an intruder within its closed neighborhood.
In this research we explore a fault-tolerant variant of identifying codes applicable to real-world systems.
Specifically, error-detecting identifying codes permit a false negative transmission from any single detector.
We investigate minimum-sized error-detecting identifying codes in several classes of graphs, including cubic graphs and infinite grids, and show that the problem of determining said minimum size in arbitrary graphs is NP-complete.
\end{abstract}

\vspace{1em}


\section{Introduction}
Let $G$ be an (undirected) graph with vertices $V(G)$ and edges $E(G)$.
The \textit{open neighborhood} of a vertex $v \in V(G)$, denoted $N(v)$, is
the set of vertices adjacent to $v$, $N(v) = \{w\in V(G): vw\in E(G)\}$.
The \textit{closed neighborhood} of a vertex $v \in V(G)$, denoted $N[v]$, is $N(v) \cup \{v\}$.
If $S \subseteq V(G)$ and every vertex in $V(G)$ is within distance 1 of some $v \in S$ (i.e., $\cup_{v \in S}{N[v]} = V(G)$), then $S$ is said to be a \emph{dominating set}; for $u \in V(G)$, we let $N_S[u] = N[u] \cap S$ and $N_S(u) = N(u) \cap S$ denote the dominators of $u$ in the closed and open neighborhoods, respectively.

A set $S \subseteq V(G)$ is called a \emph{detection system} if each vertex in $S$ is installed with a specific type of detector or sensor for locating an ``intruder" such that the set of sensor data from all detectors in $S$ can be used to precisely locate an intruder, if one is present, anywhere in the graph.
Given a detection system $S \subseteq V(G)$, two distinct vertices $u,v \in V(G)$ are said to be \emph{distinguished} if it is always possible to eliminate $u$ or $v$ as the location of an intruder (if one is present).
In order to locate an intruder anywhere in the graph, every pair of vertices must be distinguished.

Many types of detection systems with various properties have been explored throughout the years, each with their own domination and distinguishing requirements.
For example, an \emph{Identifying Code (IC)} \cite{NP-complete-ic, karpovsky} is a detection system where each detector at a vertex $v \in V(G)$ can sense an intruder within $N[v]$, but does not know the exact location.
In an IC, $S$, $u$ and $v$ are distinguished if $|N_S[u] \triangle N_S[v]| \ge 1$, where $\triangle$ denotes the symmetric difference.
A \emph{Locating-Dominating (LD) set} is a detection system that extends the capabilities of an IC by allowing a detector at vertex $v$ to differentiate an intruder in $N(v)$ from at $v$ \cite{dom-loc-acyclic, ftld}.
In an LD set, $S$, $x \in S$ is distinguished from all other vertices, and $u,v \notin S$ are distinguished if $|N_S[u] \triangle N_S[v]| \ge 1$.
Still another system is called an \emph{Open-Locating-Dominating (OLD) set}, where each detector at a vertex $v \in V(G)$ can sense an intruder within $N(v)$, but not at $v$ itself \cite{old, oldtree}.
In an OLD set, $S$, $u$ and $v$ are distinguished if $|N_S(u) \triangle N_S(v)| \ge 1$.
Lobstein \cite{dombib} maintains a bibliography of currently over 470 articles published on various types of detection systems, and other related concepts including fault-tolerant variants of ICs, LD and OLD sets.

\begin{figure}[ht]
    \centering    
    \begin{tabular}{c@{\hspace{2em}}c@{\hspace{2em}}c}
        \includegraphics[width=0.2\textwidth]{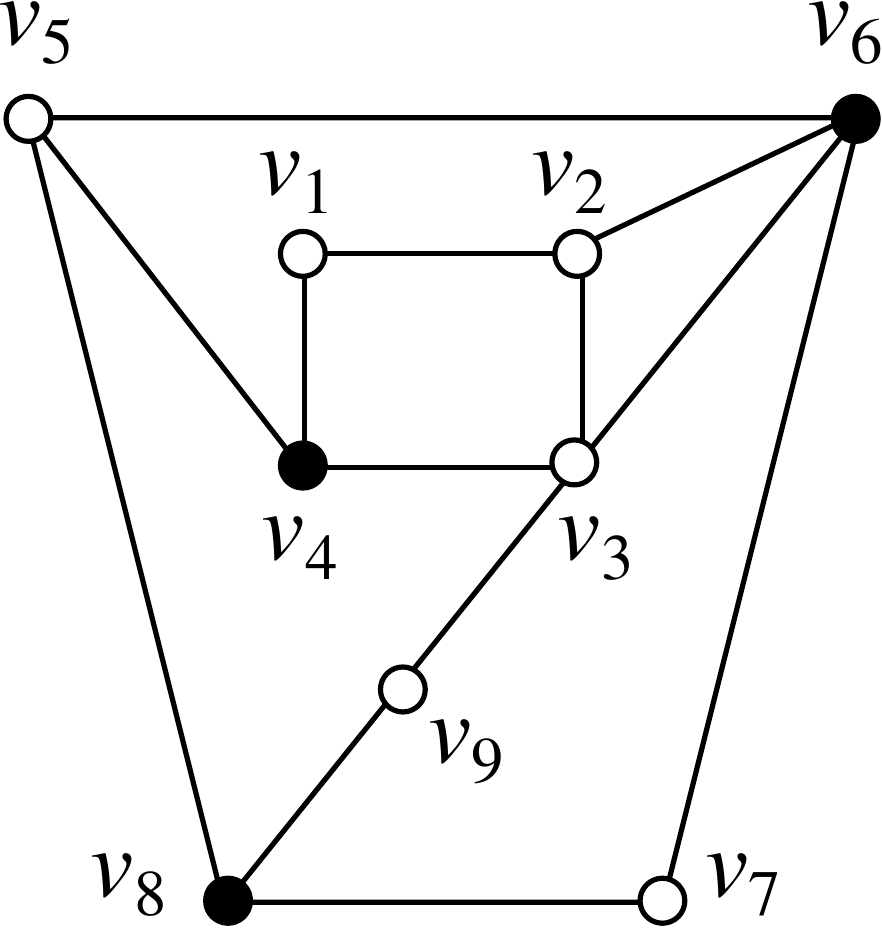} &                \includegraphics[width=0.2\textwidth]{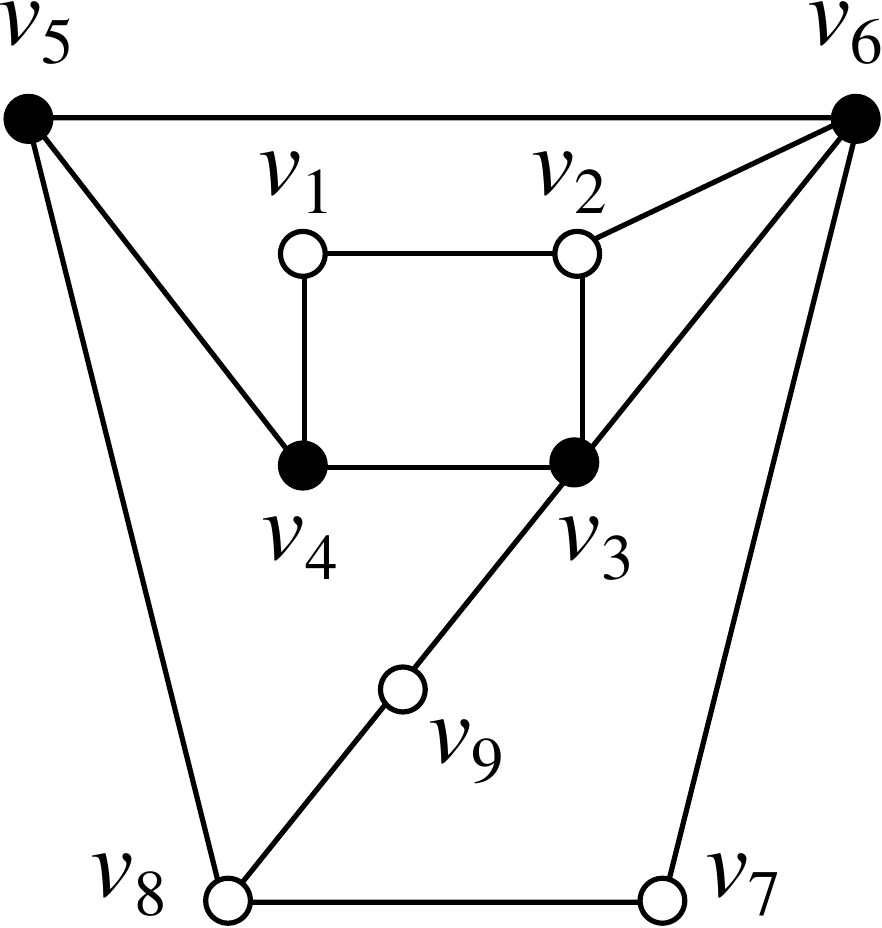} &        \includegraphics[width=0.2\textwidth]{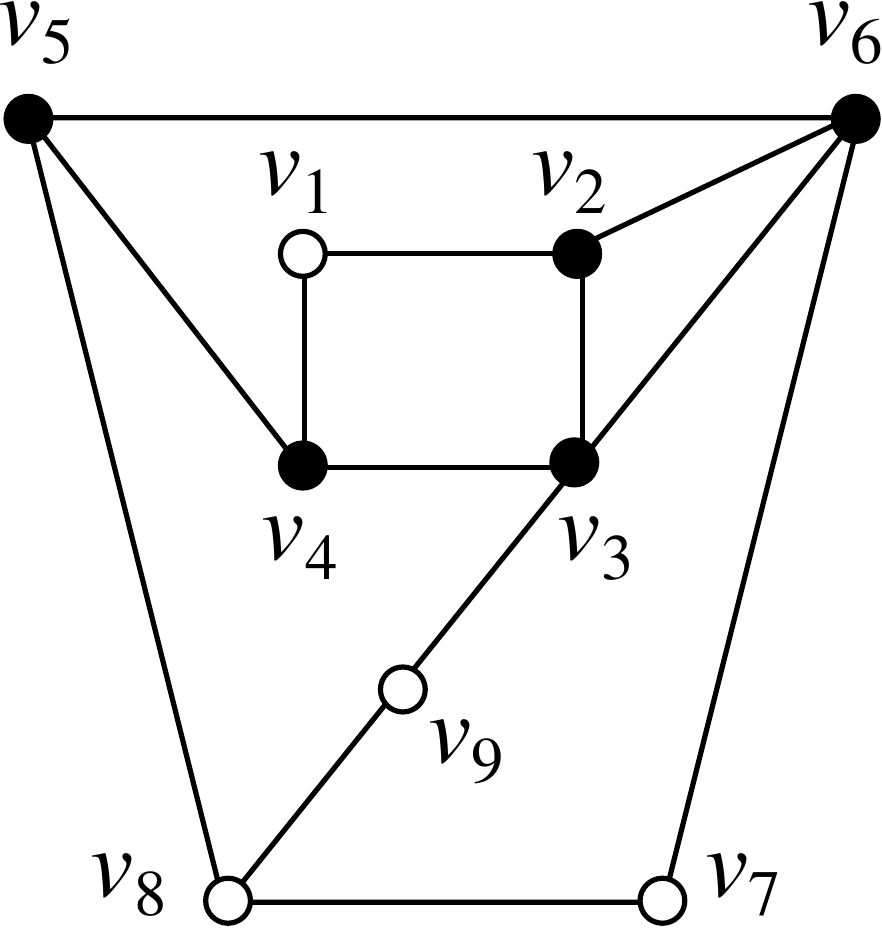} \\ (a) & (b) & (C)
    \end{tabular}
    \caption{Optimal LD (a), IC (b), and OLD (c) sets on $G_9$. Shaded vertices represent detectors.}
    \label{fig:ex-ld-ic-old}
\end{figure}

For a graph $G$, we denote the minimum cardinality of any IC, LD, or OLD set on $G$ by IC($G$), LD($G$), and OLD($G$), respectively.
Figure~\ref{fig:ex-ld-ic-old} shows LD, IC, and OLD sets on $G_9$, we can verify that there are no smaller sets for each parameter; thus, $\textrm{LD}(G_9) = 3$, $\textrm{IC}(G_9) = 4$, and $\textrm{OLD}(G_9) = 5$.

The aforementioned detection systems assume that all detectors work properly and there are no transmission errors; for applications in real-world systems, we often desire some level of fault-tolerance built into the system.
Three common fault-tolerant properties of detection systems are \emph{Redundant Detection Systems} \cite{redic, redld, ftsets}, which allow one detector to be removed, \emph{Error-Detecting Detection Systems} \cite{ourtri, detld, ftld}, which can tolerate one false negative from a sensor, and \emph{Error-Correcting Detection Systems} \cite{our3-4, errld, ft-old-cubic}, which handle any single sensor error (a false positive or false negative).

In this paper, we will focus on Error-detecting Identifying Codes (DET:ICs), including a characterization and existence criteria in Section~\ref{sec:detic-char}.
For the DET:IC parameter, DET:IC($G$) denotes the minimum cardinality of an error-detecting IC on graph $G$.
For many detection systems and their fault tolerant variants, minimizing a detection system is known to be NP-complete for arbitrary graphs \cite{NP-complete-ic, NP-complete-ld, errld, redld, detld, redic, old}.
In Section~\ref{sec:npc}, we will prove the problem of determining DET:IC(G) for an arbitrary graph $G$ is also NP-complete.
In Section \ref{sec:special} we investigate minimum-sized error-detecting identifying codes in several classes of graphs, including cubic graphs and infinite grids

\section{Properties of DET:IC}\label{sec:detic-char}

Detection systems commonly use general terminology such as ``dominated" or ``distinguished", whose specific definitions vary depending on the sensors' capabilities and the level of fault-tolerance.
The following definitions are specifically for identifying codes and their fault-tolerant variants; assume that $S \subseteq V(G)$ is the set of detectors.
Further, existence arguments assume $S = V(G)$, as extraneous sensors cannot remove information from the system.

\begin{definition}\label{def:k-dom}
A vertex $v \in V(G)$ is \emph{$k$-dominated} by a dominating set $S$ if $|N_S[v]| = k$.
\end{definition}

\begin{definition}\label{def:k-disty}
If $S$ is a dominating set and $u,v \in V(G)$, $u$ and $v$ are \emph{$k$-distinguished} if $|N_S[u] \triangle N_S[v]| \ge k$, where $\triangle$ denotes the symmetric difference.
\end{definition}

\begin{definition}\label{def:k-sharp-disty}
If $S$ is a dominating set and $u,v \in V(G)$, $u$ and $v$ are \emph{$k^\#$-distinguished} if $|N_S[u] - N_S[v]| \ge k$ or $|N_S[v] - N_S[u]| \ge k$.
\end{definition}

We will also use terms such as ``at least $k$-dominated" to denote $j$-dominated for some $j \ge k$.

\vspace{0.6em}
Jean and Seo \cite{redic} have shown the necessary and sufficient properties of a fault-tolerant identifying code: redundant identifying code (RED:IC).
Seo and Slater \cite{separating} characterized separating sets, which are more general, set-theoretic forms of fault-tolerant detection systems; we can convert their characterization to the following for error-detecting identifying codes (DET:ICs).

\begin{definition}\label{def:ic-char}
A detector set, $S \subseteq V(G)$, is an IC if and only if each vertex is at least 1-dominated and all pairs are 1-distinguished.
\end{definition}

\begin{theorem}[\cite{separating}]\label{theo:detic-char}
A detector set, $S \subseteq V(G)$, is a DET:IC if and only if each vertex is at least 2-dominated and all pairs are $2^\#$-distinguished.
\end{theorem}

For finite graphs, we use the notations $\textrm{IC}(G)$, and $\textrm{DET:IC}(G)$ to denote the cardinality of the smallest possible such sets on graph $G$, respectively.
For infinite graphs, we measure via the \emph{density} of the subset, which is defined as the ratio of the size of the subset to the size of the whole set \cite{ourtri,ftsets}.
Formally, for locally-finite (i.e. $B_r(v)$ finite for finite $r$) $G$, this is defined as $\limsup_{r \to \infty}{\frac{|B_r(v) \cap S|}{|B_r(v)|}}$ for any $v \in V(G)$, where $B_r(v) = \{ u \in V(G) : d(u,v) \le r\}$ denotes the ball with radius $r$ around $v$.
We use the notations $\textrm{IC\%}(G)$, and $\textrm{DET:IC\%}(G)$ to denote the lowest density of any possible such set on $G$ \cite{ourtri,ftsets}.
Note that density is also defined for finite graphs.

\begin{figure}[ht]
    \centering    
    \begin{tabular}{c@{\hspace{4em}}c}
        \includegraphics[width=0.2\textwidth]{fig/ex-ic.pdf} &        \includegraphics[width=0.2\textwidth]{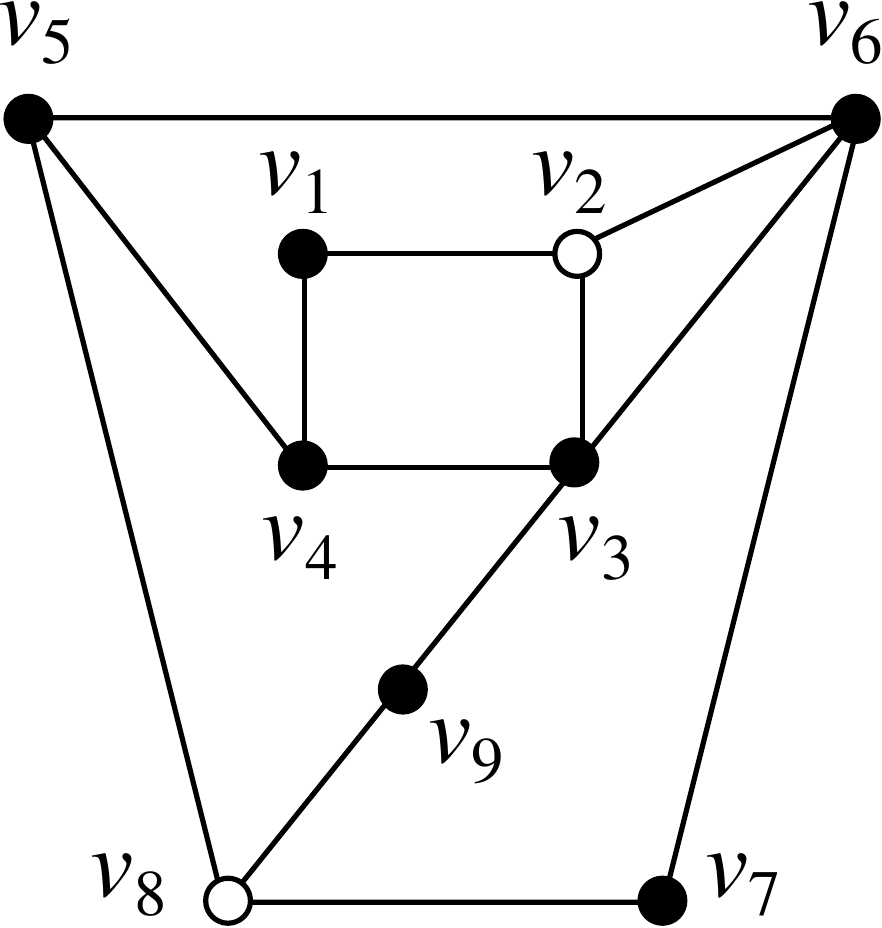} \\ (a) & (b) 
    \end{tabular}
    \caption{Optimal IC (a) and DET:IC (b) sets on $G_9$. Shaded vertices represent detectors.}
    \label{fig:ex-ic-detic}
\end{figure}

Figure~\ref{fig:ex-ic-detic} shows an example of IC and DET:IC on $G_9$.
In the IC set (a), we see that every vertex is at least 1-dominated.
We see that $N_S[v_1] \triangle N_S[v_2] = \{ v_3,v_4,v_6\}$, $N_S[v_1] \triangle N_S[v_3] = \{v_3,v_6\}$, $N_S[v_2] \triangle N_S[v_6] = \{v_5\}$, and so on; all vertex pairs are 1-distinguished.
Therefore, Definition~\ref{def:ic-char} yields that (a) is an IC.
For the DET:IC set in (b) clearly all vertices are at least 2-dominated.
We also see all pairs are $2^\#$-distinguished.
For example, $(v_1, v_2)$ are $2^\#$-distinguished because $N_S[v_2] - N_S[v_1] = \{v_3,v_6\}$ and $(v_2, v_4)$ are $2^\#$-distinguished because $N_S[v_4] - N_S[v_2] = \{v_4,v_5\}$.
Thus, by Theorem~\ref{theo:detic-char} it is a DET:IC.
There are no smaller IC or DET:IC sets exist on $G_9$; thus, $\textrm{IC}(G_9) = 4$ and $\textrm{DET:IC}(G_9) = 7$.
If we would prefer to use densities, we also have that $\textrm{IC\%}(G_9) = \frac{4}{9}$ and $\textrm{DET:IC\%}(G_9) = \frac{7}{9}$.

\begin{observation}\label{obs:detic-no-isolated}
If $G$ has a DET:IC, then $G$ lacks isolated vertices.
\end{observation}

\begin{definition}\cite{ld-twin-free} 
Two distinct vertices $u,v \in V(G)$ are said to be \emph{twins} if $N[u] = N[v]$ (\emph{closed twins}) or $N(u) = N(v)$ (\emph{open twins}).
\end{definition}

It is easy to see $G$ has an IC if and only if $G$ has no closed-twins.
If two distinct vertices $u,v \in V(G)$ are open-twins, then $|N[u] - N[v]| = 1$ and $|N[v] - N[u]| = 1$; thus it is impossible to $2^\#$-distinguish $u$ and $v$.

\begin{observation}\label{obs-detic-twinfree}
If $G$ has a DET:IC, then $G$ is twin-free.
\end{observation}

\begin{lemma}\label{lemma:triangle-dom}
If $G$ has a DET:IC, then for each $uv \in E(G)$, $deg(u) \ge t + 3$ or $deg(v) \ge t + 3$, where $t$ is the number of triangles in $G$ containing the edge $uv$.
\end{lemma}
\begin{proof}
Suppose for a contradiction that $G$ has a DET:IC but $uv \in E(G)$ with $deg(u) \le t + 2$ and $deg(v) \le t + 2$.
Consider vertex $u$: the presence of $t$ triangles and the $uv$ edge imply that $t + 1 \le deg(u) \le t + 2$; however these $t+1$ vertices in $\{v\} \cup (N(u) \cap N(v))$ and $u$ itself do not contribute to distinguishing $u$ and $v$ since they are contained in $N[u] \cap N[v]$.
Thus, $|N[u] - N[v]| \le (t + 3) - (t + 2) = 1$, and by symmetry $|N[v] - N[u]| \le 1$, contradicting that $G$ has a DET:IC.
\end{proof}

\begin{theorem}\label{theo:detic-exist-alt}
A graph $G$ has a DET:IC if and only if it satisfies the following properties.
\begin{enumerate}[label=\roman*,noitemsep]
    \item $G$ has $\delta(G) \ge 1$
    \item $G$ is open-twin-free
    \item each $uv \in E(G)$ has $deg(u) \ge t + 3$ or $deg(v) \ge t + 3$, where $t$ is the number of triangles containing the edge $uv$.
\end{enumerate}
\end{theorem}
\begin{proof}
If $G$ has a DET:IC, then Observations \ref{obs:detic-no-isolated} and  \ref{obs-detic-twinfree}, and Lemma~\ref{lemma:triangle-dom} give us the required properties; thus, we need only prove the converse.
Specifically, we will show that if $G$ satisfies the above properties, then $S = V(G)$ is a DET:IC for $G$.
By the requirement of $\delta(G) \ge 1$, we know that every vertex is at least 2-dominated.
We now need only show that each vertex pair is distinguished; let $u,v \in V(G)$ be distinct vertices.
If $d(u,v) \ge 3$, then $u$ and $v$ cannot share common dominators and are each at least $2$-dominated; thus, $u$ and $v$ are distinguished.
If $d(u,v) = 2$, then the fact that $G$ is open-twin-free lets us assume by symmetry that $\exists p \in N(u) - N[v]$, meaning $u$ and $v$ are distinguished by $u$ and $p$.
Finally, we assume that $d(u,v) = 1$, which implies that $uv \in E(G)$.
By property~\toroman{3}, we know that $deg(u) \ge t + 3$ or $deg(v) \ge t + 3$; without loss of generality, assume $deg(u) \ge t + 3$.
Because $uv \in E(G)$, there are precisely $t$ vertices in $N(u) \cap N(v)$; thus, $deg(u) \ge t + 3$ implies that $|N[u] - N[v]| \ge 2$, so $u$ and $v$ are distinguished, completing the proof.
\end{proof}

\begin{corollary}\label{cor:detic-Cn}
$C_n$ do not permit DET:IC for any $n$.
\end{corollary}

If $G$ is a cubic graph and $uv \in E(G)$, then Theorem~\ref{theo:detic-exist-alt} requires $deg(u) \ge t + 3$ or $deg(v) \ge t + 3$.
However, $G$ is cubic, so it must be that $t = 0$.
Because $uv \in E(G)$ was selected arbitrarily, this implies that $G$ is triangle-free.

\begin{corollary}\label{cor:detic-exist-alt-cubic}
A cubic graph $G$ has a DET:IC if and only if $G$ is open-twin-free and triangle-free.
\end{corollary}

If two distinct vertices $u,v \in V(G)$ are leaf vertices with the common support vertex, then they are open-twins; thus it is impossible to $2^\#$-distinguish $u$ and $v$.

\begin{observation}\label{obs-detic-support}
If $G$ has a DET:IC, then every support vertex in $G$ is associated with a unique leaf child.
\end{observation}

\begin{theorem}\label{theo:detic-tree}
There are no trees with DET:IC.
\end{theorem}
\begin{proof}
Suppose there is a tree $T$ with DET:IC.
Consider the sub-graph $T' \subseteq T$ which excludes all leaves from $T$.
Note that $|T'| \ge 1$ because the only tree consisting entirely of leaves is $P_2$, which does not permit DET:IC.
Let $x$ be a support vertex in $T$ (or equivalently a leaf vertex in $T'$), then by Lemma~\ref{lemma:triangle-dom} we have $deg(x) \ge 3$ in $T$.
By Observation~\ref{obs-detic-support}, each support vertex has exactly one leaf, implying $deg(x) \ge 2$ in $T'$, contradicting that $x$ is a leaf in $T'$.
\end{proof}

Suppose $G$ has a DET:IC and a 4-cycle $abcd$.
Applying Lemma~\ref{lemma:triangle-dom} to the edges of $abcd$ gives by symmetry that $deg(a) \ge 3$ and $deg(c) \ge 3$.
However, if $deg(b) = deg(d) = 2$, then $b$ and $d$ would be open twins, so by symmetry $deg(b) \ge 3$ is required as well.

\begin{observation}\label{obs-detic-c4}
If $G$ has a DET:IC, then each 4-cycle in $G$ can have at most one vertex with degree 2.
\end{observation}

\begin{figure}[ht]
    \centering
    \includegraphics[width=0.15\textwidth]{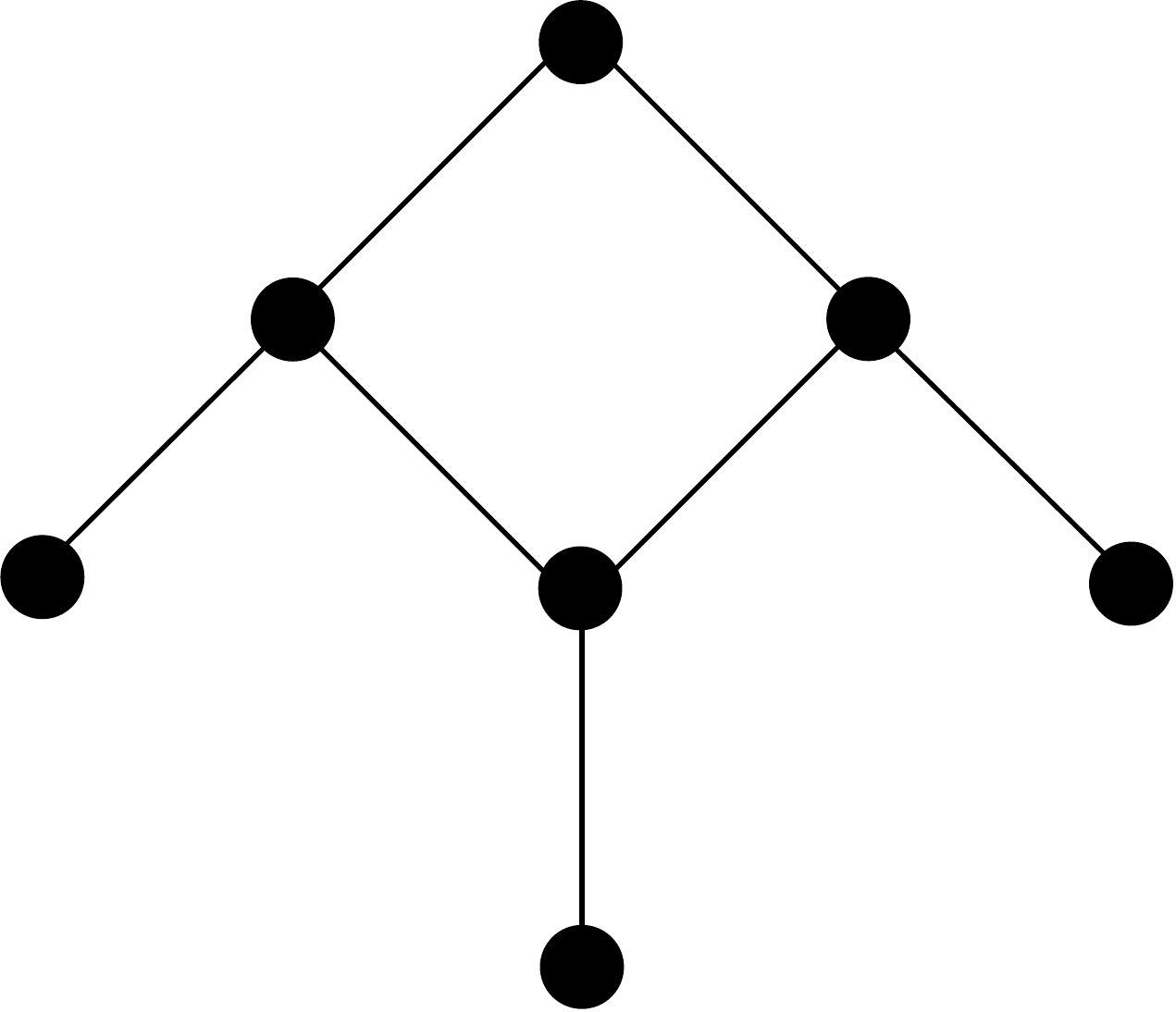}
    \caption{$G_{7,7}$ with $DET:IC(G_{7,7}) = 7$}
    \label{fig:g7}
\end{figure}
\FloatBarrier

\begin{theorem}\label{theo:detic-g7-7-unique}
Let $G_{n,m}$ have $n$ vertices and $m$ edges.
Then $G_{7,7}$, as shown in Figure~\ref{fig:g7}, is the first graph with DET:IC in the lexicographic ordering of $(n,m)$ tuples; i.e., the graph with the smallest number of edges given the smallest number of vertices.
\end{theorem}
\begin{proof}
Let $S=V(G)$ be a DET:IC.
By Theorem~\ref{theo:detic-tree}, $G$ contains at least one cycle and we consider the existence of cycles of a given size.

Case 1: $G$ has a 3-cycle $abc$.
By applying Lemma~\ref{lemma:triangle-dom} to the edges of $abc$, by symmetry we can assume that $deg(a) \ge 4$ and $deg(b) \ge 4$, so let $a',a'' \in N(a) - \{b,c\}$ and $b',b'' \in N(b) - \{a,c\}$ with $a' \neq a''$ and $b' \neq b''$.
If $a' = b'$ and $a'' = b''$, then edge $ab \in E(G)$ has $t \ge 3$, so Lemma~\ref{lemma:triangle-dom} gives that, without loss of generality, $deg(a) \ge 6$, implying $n \ge 7$ but $m > 7$, and so we are done.
Otherwise, we can assume by symmetry that $a'' \neq b''$.
If $a' = b'$ then applying Lemma~\ref{lemma:triangle-dom} to $ab \in E(G)$ yields $t \ge 2$, implying without loss of generality that $deg(a) \ge 5$, so we again have $n \ge 7$ but $m > 7$ and so are done.
Otherwise, we can assume that $a' \neq b'$, so $\{a',a'',b',b''\}$ are all distinct; then $n = 7$ and $m = 7$, but $(a',a'')$ and $(b',b'')$ are twins, contradicting that DET:IC exists.

Case 2: $G$ has a 4-cycle $abcd$. We can assume $G$ is triangle-free, as otherwise we fall into case 1 and would be done.
By applying Observation~\ref{obs-detic-c4} to $abcd$, we can assume that $a' \in N(a) - \{b,d\}$, $b' \in N(b) - \{a,c\}$, and $c' \in N(c) - \{b,d\}$.
Because $G$ is now assumed to be triangle-free, we know that $a' \neq b'$ and $b' \neq c'$.
If $a' = c'$ then we have $6$ vertices, but $(a,c)$ currently are open twins, so $n \ge 7$ will be required to distinguish them, but $m > 7$, so we are done.
Otherwise, we can assume that $a' \neq c'$, so $\{a',b',c'\}$ are distinct.
We now have $n=7$ and $m=7$, so no more vertices or edges may be added, and we see that this is exactly the $G_{7,7}$ graph shown in Figure~\ref{fig:g7}.

Case 3: $G$ has a cycle, $C_k$, of length $5 \le k \le 7$, and $G$ has girth $k$.
By applying Lemma~\ref{lemma:triangle-dom} to the edges of $C_k$, we have that there must be at least $\ell = \ceil{\frac{k}{2}}$ vertices in $C$ which have degree at least 3.
If any of these $\ell$ additional vertices are not distinct, we would contradict that $G$ has girth $k$, so we can assume these $\ell$ additional vertices are distinct.
Then $n \ge k + \ell > 7$ because $k \ge 5$, a contradiction, completing the proof.
\end{proof}

\begin{observation}\label{obs:det-ic-upper}
For each $n \ge 7$, there exists a graph $G_n$ on $n$ vertices which has $\textrm{DET:IC}(G_n) = n$.
\end{observation}

\begin{figure}[ht]
    \centering
    \begin{tabular}{c@{\hspace{1em}}c@{\hspace{1em}}c}
        \centered{\includegraphics[width=0.15\textwidth]{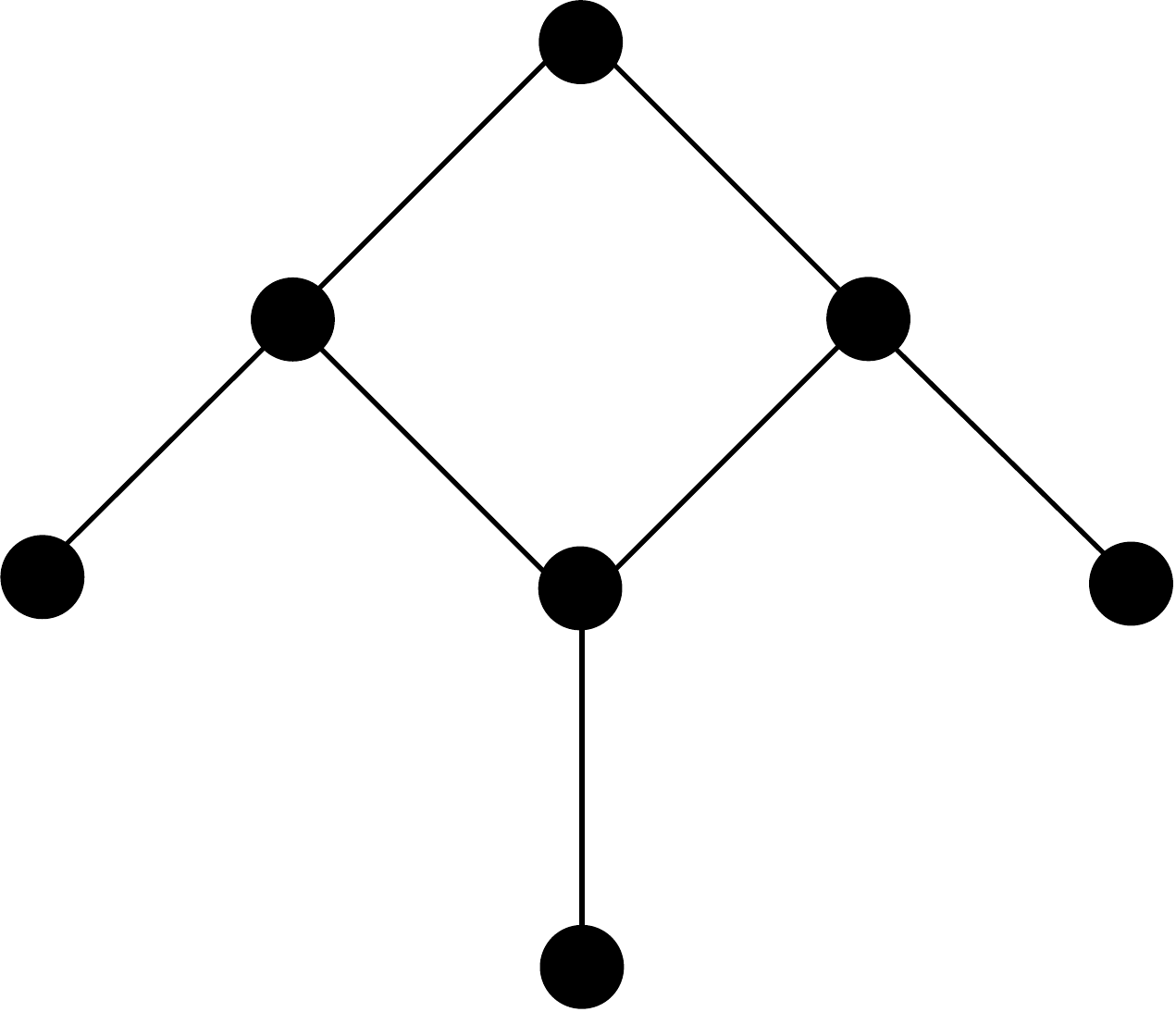}} &
        \centered{\includegraphics[width=0.15\textwidth]{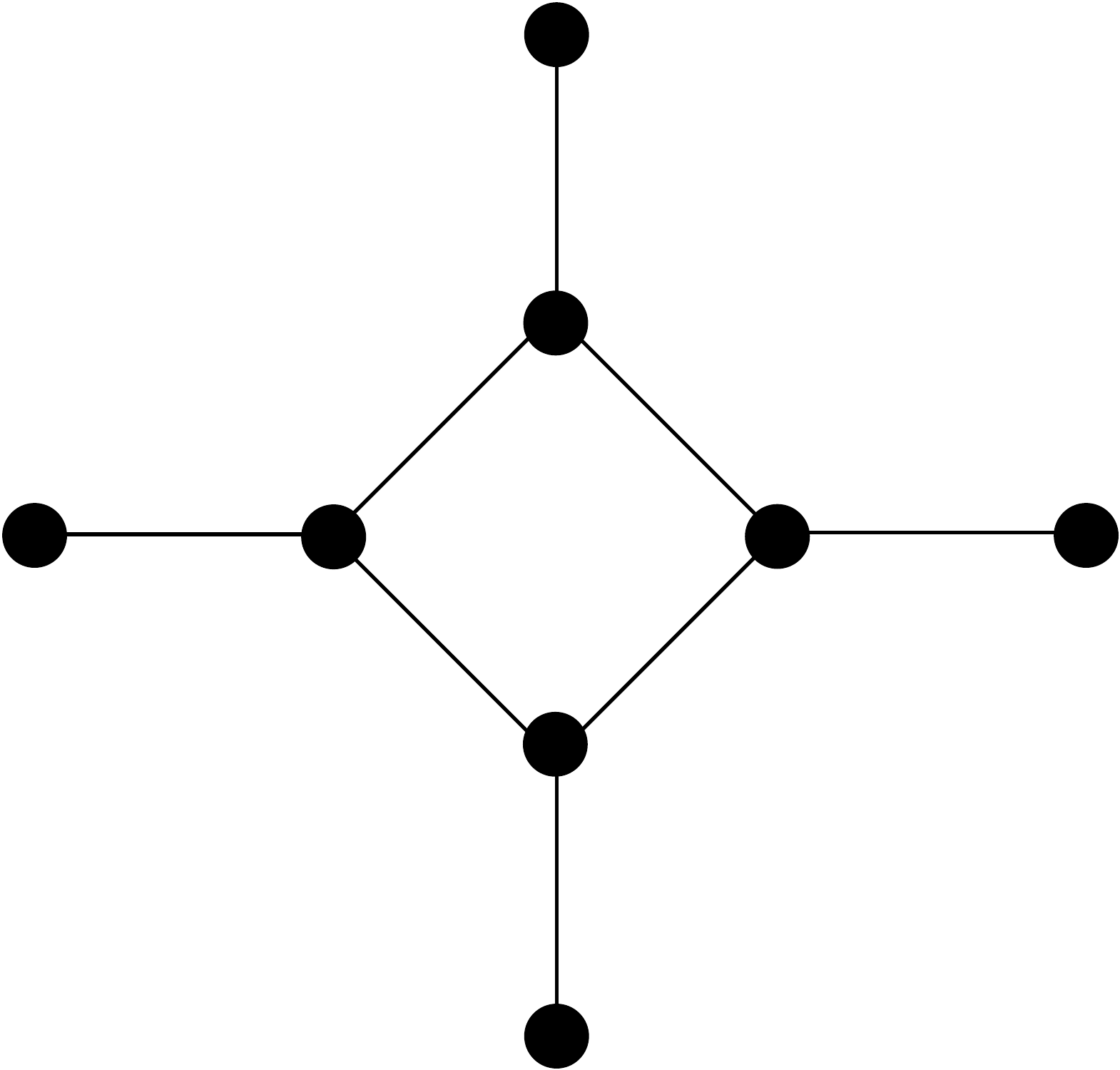}} &
        \centered{\includegraphics[width=0.12\textwidth]{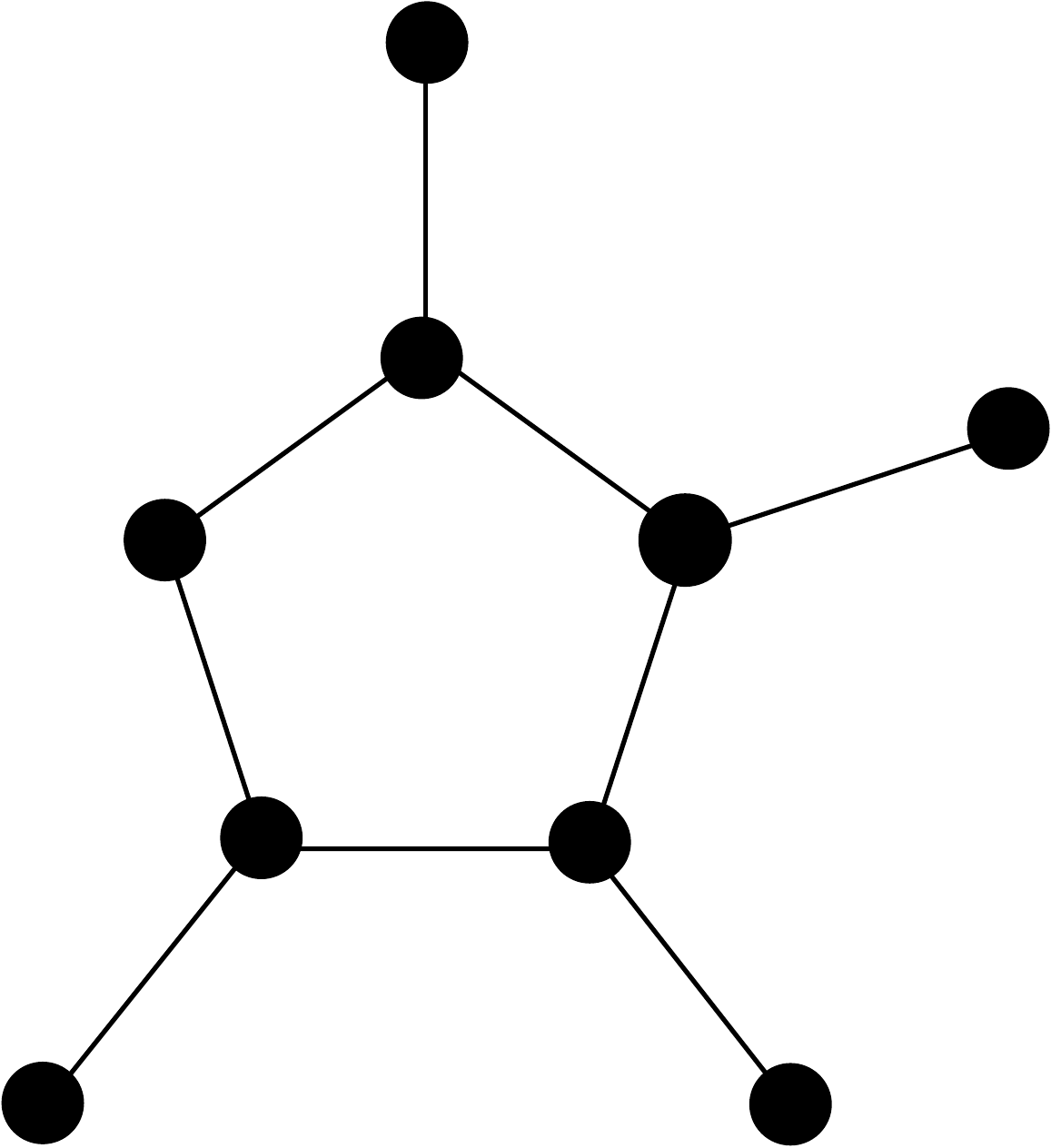}} \\\\
    \end{tabular}
    \begin{tabular}{c@{\hspace{1em}}c@{\hspace{1em}}c}
        \centered{\includegraphics[width=0.15\textwidth]{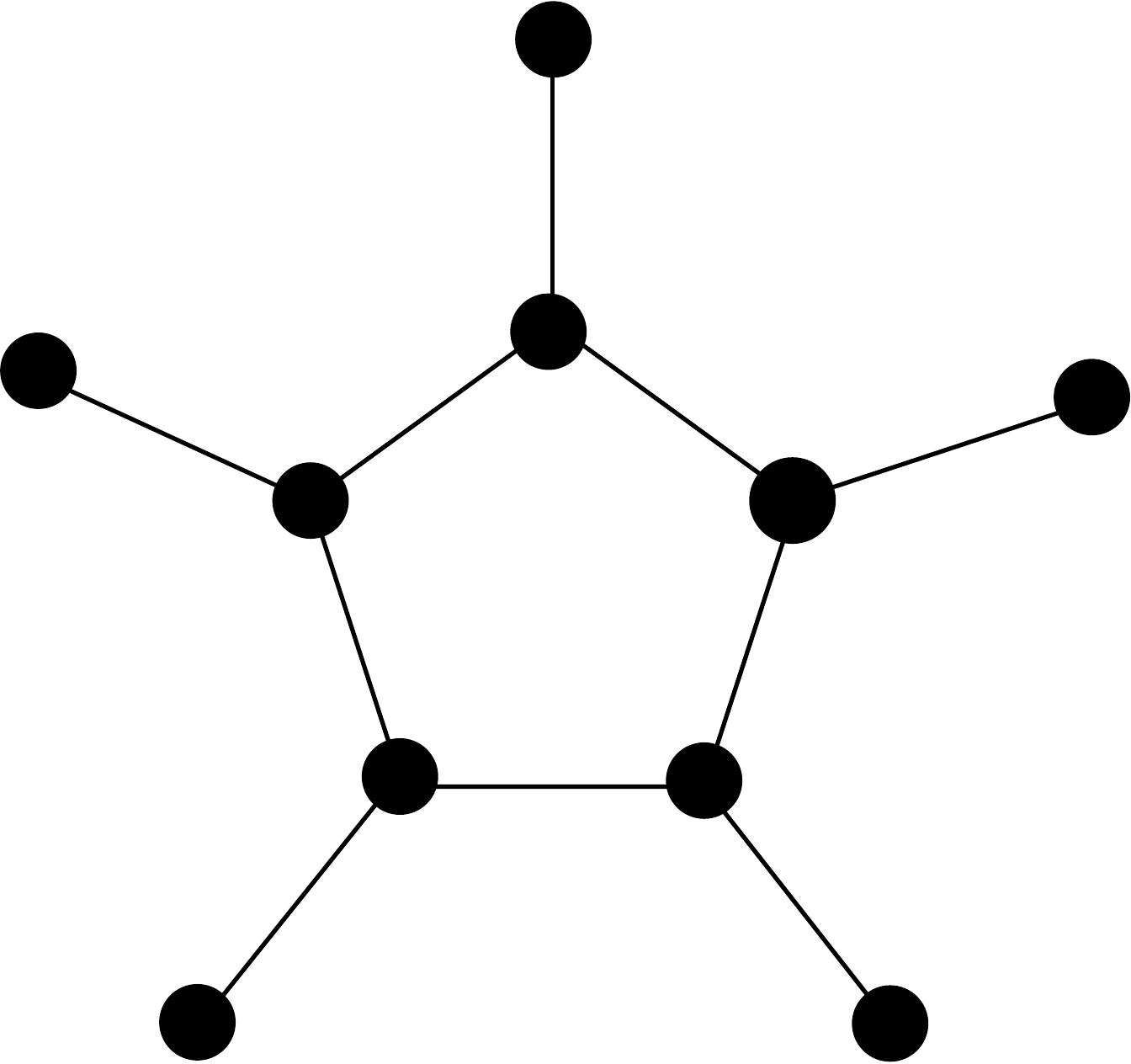}} &
        \centered{\includegraphics[width=0.12\textwidth]{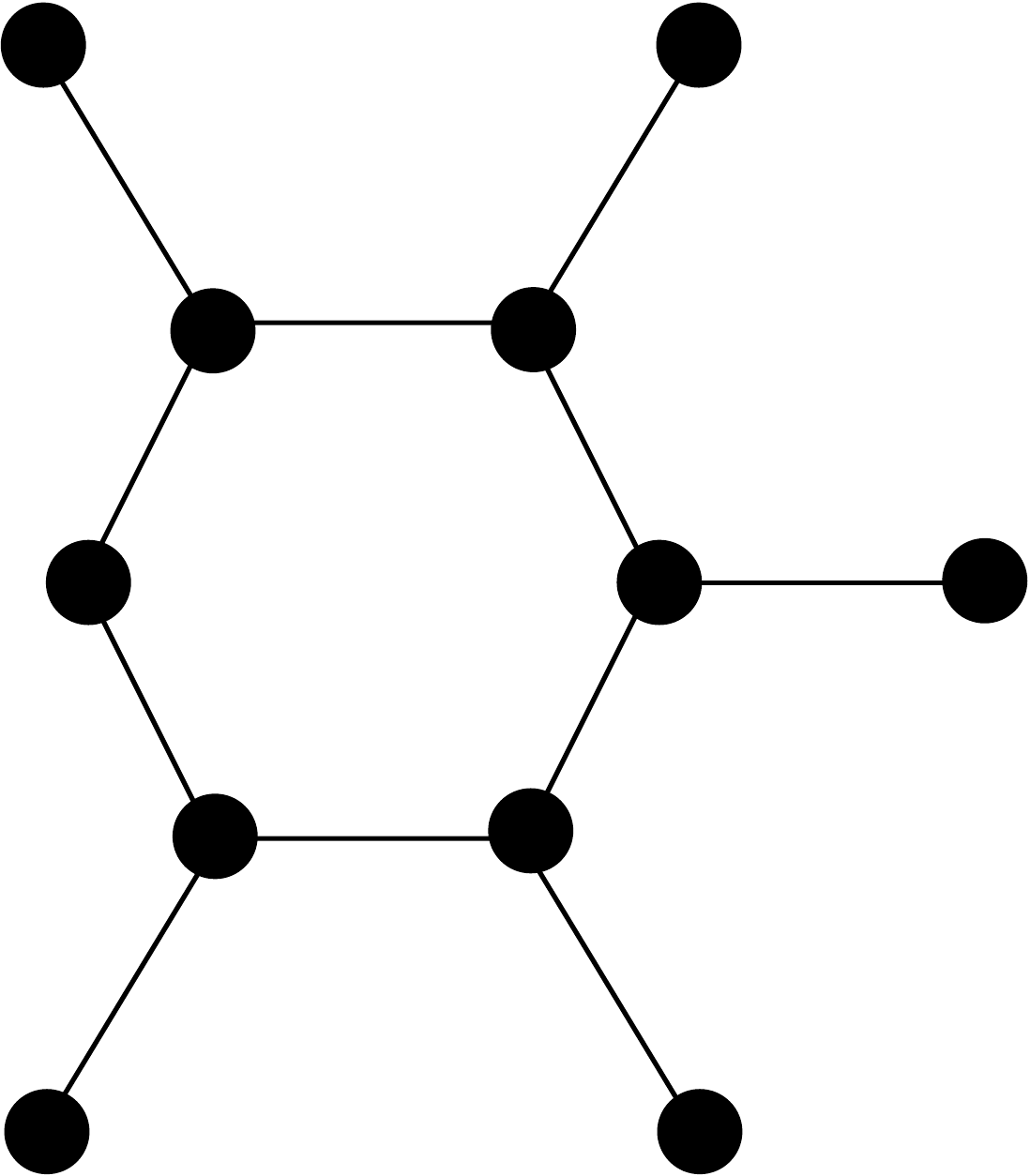}} &
        \centered{\includegraphics[width=0.15\textwidth]{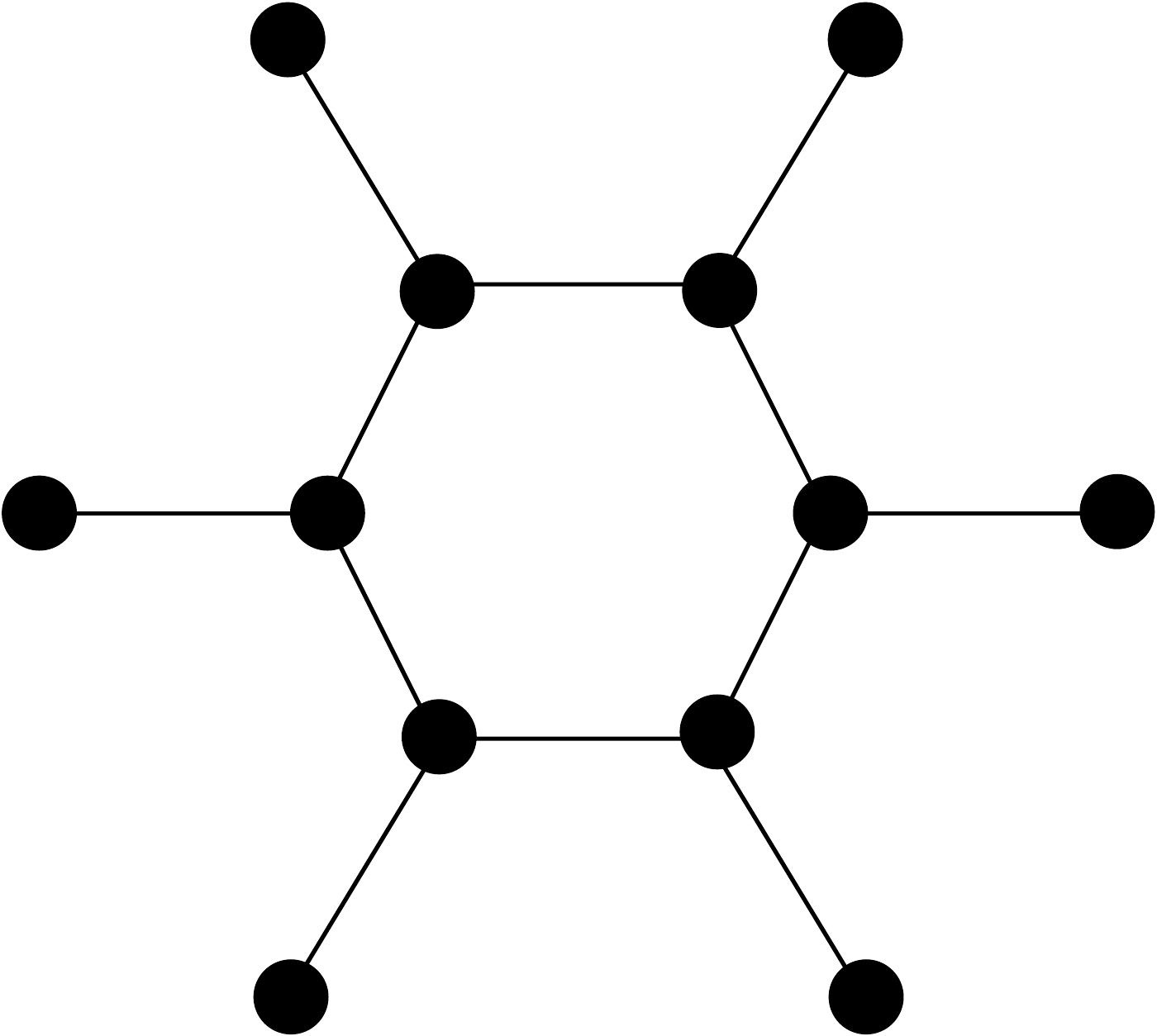}} 
    \end{tabular}
    \caption{A family of graphs with $\textrm{DET:IC}(G)=n$.}
    \label{fig:det-ic-upper-fam}
\end{figure}

Note that Observation~\ref{obs:det-ic-upper} gives an infinite family of graphs which have $\textrm{DET:IC}(G) = n$.
Figure~\ref{fig:det-ic-upper-fam} gives example from this family for $7 \le n \le 12$.

\FloatBarrier

\section{NP-completeness of Error-detecting IC}\label{sec:npc}

It has been shown that many graphical parameters related to detection systems, such as finding the cardinality of the smallest IC, LD, or OLD sets, are NP-complete problems \cite{ld-ic-np-complete-2, NP-complete-ic, NP-complete-ld, old}.
We will now prove that the problem of determining the smallest DET:IC set is also NP-complete.
For additional information about NP-completeness, see Garey and Johnson \cite{np-complete-bible}.

\npcompleteproblem{3-SAT}{Let $X$ be a set of $N$ variables.
Let $\psi$ be a conjunction of $M$ clauses, where each clause is a disjunction of three literals from distinct variables of $X$.}{Is there is an assignment of values to $X$ such that $\psi$ is true?}

\npcompleteproblem{Error-Detecting Identifying Code (DET-IC)}{A graph $G$ and integer $K$ with $2????? \le K \le |V(G)|$.}{Is there a DET:IC set $S$ with $|S| \le K$? Or equivalently, is DET:IC($G$) $\le K$?}

\begin{theorem}
The DET-IC problem is NP-complete.
\end{theorem}

\begin{wrapfigure}{r}{0.45\textwidth}
    \centering
    \includegraphics[width=0.45\textwidth]{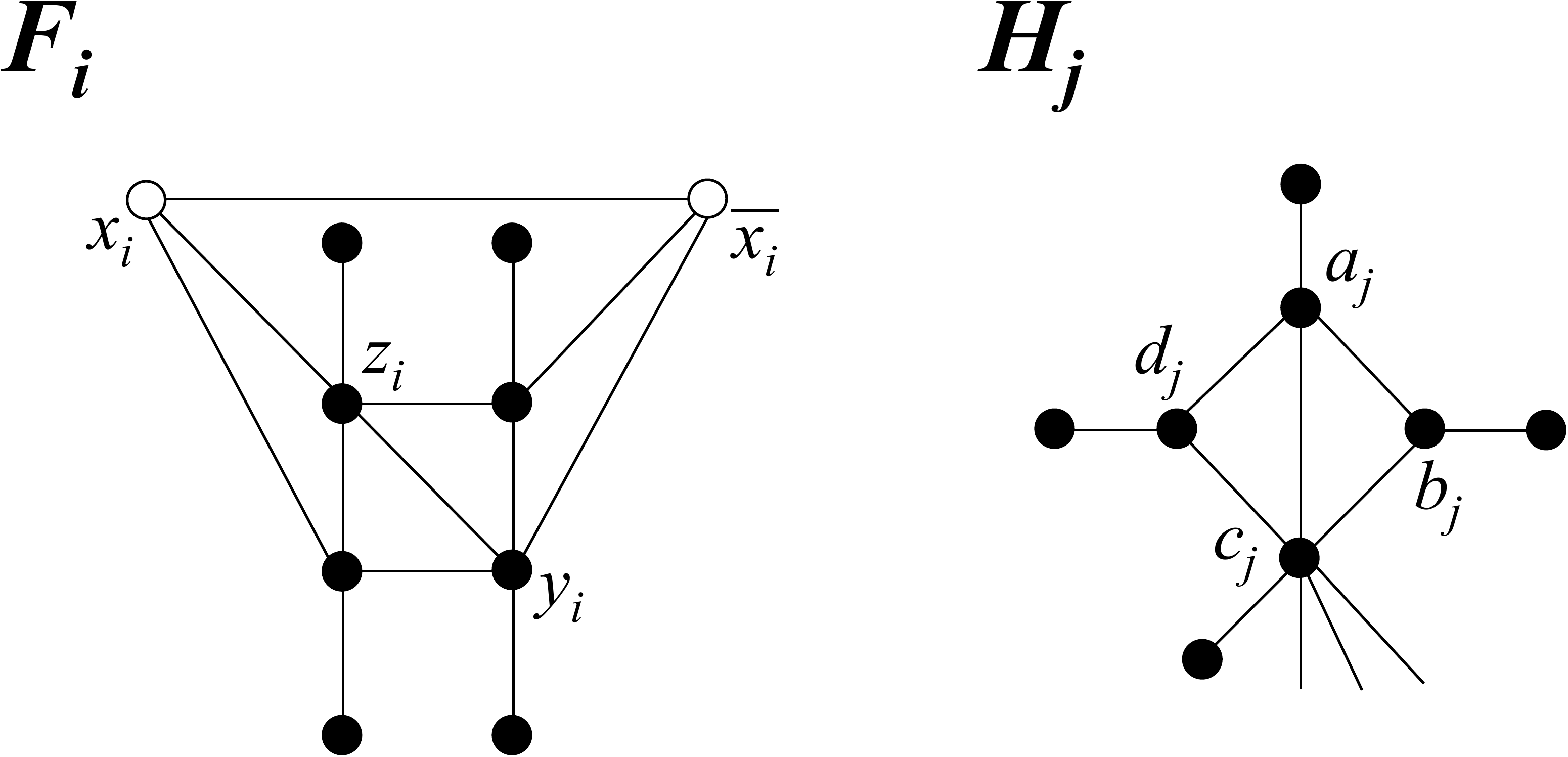}
    \caption{Variable and Clause graphs}
    \label{fig:variable-clause}
\end{wrapfigure}
\cbeginproof
Clearly, DET-IC is NP, as every possible candidate solution can be generated nondeterministically in polynomial time (specifically, $O(n)$ time), and each candidate can be verified in polynomial time using Theorem~\ref{theo:detic-char}.
To complete the proof, we will now show a reduction from 3-SAT to DET-IC.

\begin{figure}[ht]
    \centering
    \includegraphics[width=0.9\textwidth]{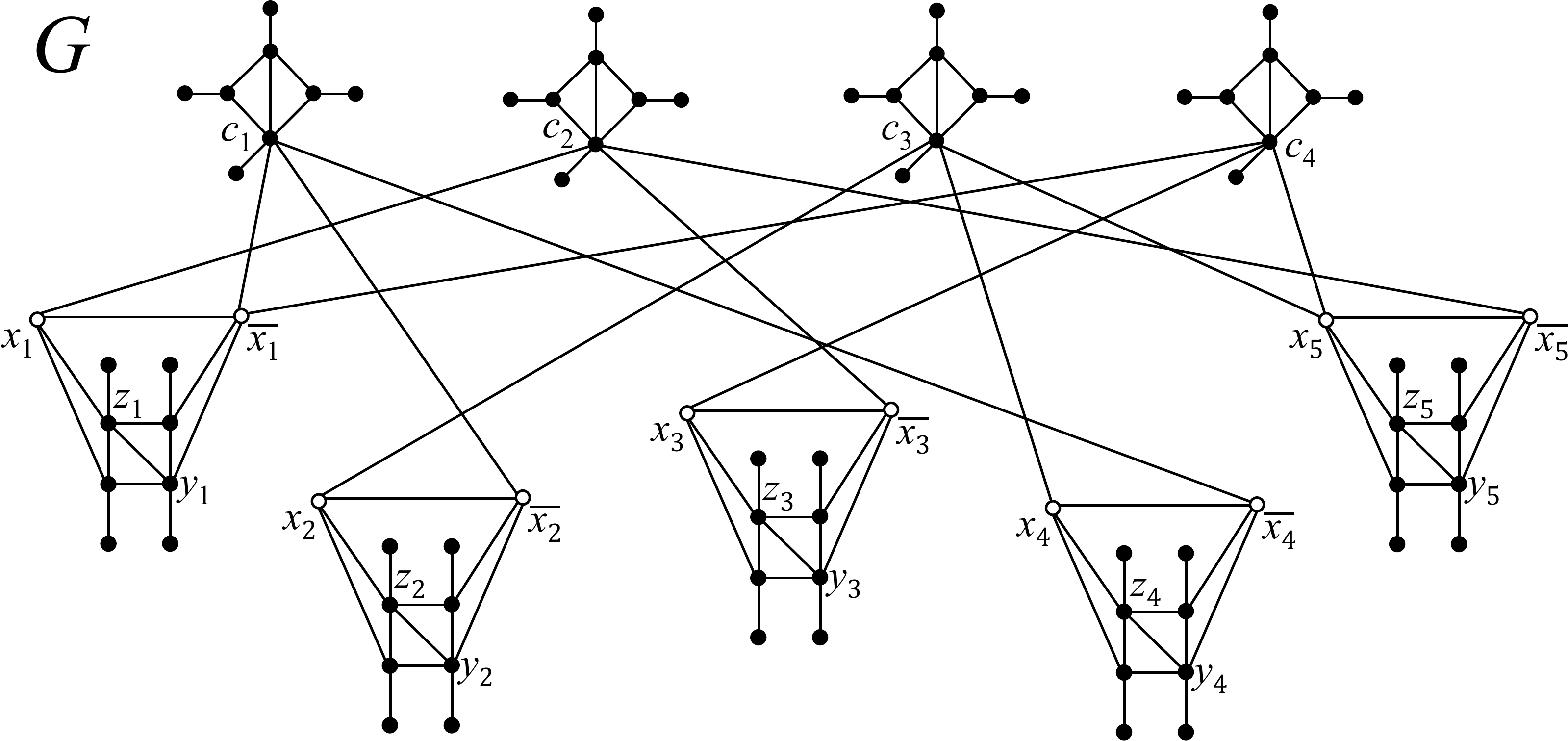}
    \caption{\begin{tabular}[t]{l} Construction of $G$ from $(\overline{x_1} \lor \overline{x_2} \lor \overline{x_4}) \land (x_1 \lor \overline{x_3} \lor \overline{x_5}) \land (x_2 \lor x_4 \lor x_5) \land (\overline{x_1} \lor x_3 \lor x_5)$ \protect\\ with $N = 5$, $M = 4$, $K = 47$ \end{tabular}}
    \label{fig:example-clause}
\end{figure}

Let $\psi$ be an instance of the 3-SAT problem with $M$ clauses on $N$ variables.
We will construct a graph, $G$, as follows.
For each variable $x_i$, create an instance of the $F_i$ graph (Figure~\ref{fig:variable-clause}); this includes a vertex for $x_i$ and its negation $\overline{x_i}$.
For each clause $c_j$ of $\psi$, create a new instance of the $H_j$ graph (Figure~\ref{fig:variable-clause}).
For each clause $c_j = \alpha \lor \beta \lor \gamma$, create an edge from the $c_j$ vertex to $\alpha$, $\beta$, and $\gamma$ from the variable graphs, each of which is either some $x_i$ or $\overline{x_i}$; for an example, see Figure~\ref{fig:example-clause}.
The resulting graph has precisely $10N + 8M$ vertices and $14N + 12M$ edges, and can be constructed in polynomial time.
To complete the problem instance, we define $K = 9N + 8M$.

Suppose $S \subseteq V(G)$ is a DET:IC on $G$ with $|S| \le K$.
By Theorem~\ref{theo:detic-char}, every vertex must be at least 2-dominated; thus, we require at least $8N + 8M$ detectors, as shown by the shaded vertices in Figure~\ref{fig:variable-clause}.
Additionally, in each $F_i$ we see that $y_i$ and $z_i$ are not distinguished unless $\{x_i,\overline{x_i}\} \cap S \neq \varnothing$.
Thus, we find that $|S| \ge 9N + 8M = K$, implying that $|S| = K$, so $|\{x_i,\overline{x_i}\} \cap S| = 1$ for each $i \in \{1, \hdots, N \}$.
For each $H_j$, we see that $a_j$ and $c_j$ are not distinguished unless $c_j$ is adjacent to at least one additional detector vertex.
As no more detectors may be added, it must be that each $c_j$ is now dominated by one of its three neighbors in the $F_i$ graphs; therefore, $\psi$ is satisfiable.

For the converse, suppose we have a solution to the 3-SAT problem $\psi$; we will show that there is a DET:IC, $S$, on $G$ with $|S| \le K$.
We construct $S$ by first including all of the $8N + 8M$ vertices needed for 2-domination.
Then, for each variable, $x_i$, if $x_i$ is true then we let the vertex $x_i \in S$; otherwise, we let $\overline{x_i} \in S$.
Thus, the fully-constructed $S$ has $|S| = 9N + 8M = K$.
Because we selected each $x_i \in S$ or $\overline{x_i} \in S$ based on a satisfying truth assignment for $\psi$, each $c_j$ must be adjacent to at least one additional detector vertex from the $F_i$ graphs.
It can then be shown that all vertex pairs are distinguished, so $S$ is a DET:IC for $G$ with $|S| \le K$, completing the proof.
\cendproof

\FloatBarrier

\section{DET:IC in Special classes of graphs}\label{sec:special}

\subsection{Hypercubes}

\begin{figure}[ht]
    \centering
    \includegraphics[width=0.45\textwidth]{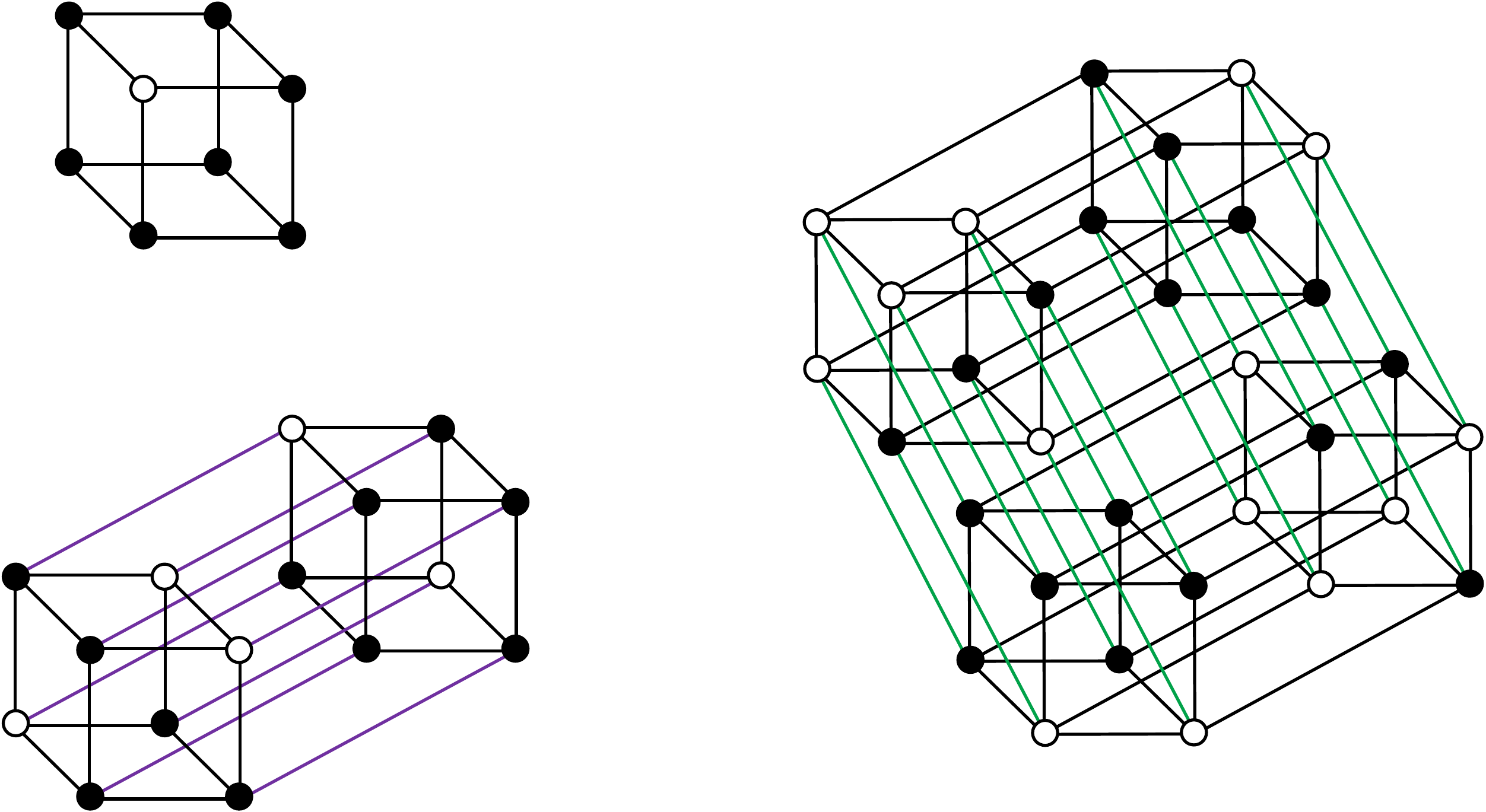}
    \caption{DET:ICs for $Q_n$ with $n \le 5$}
    \label{fig:hypercubes}
\end{figure}

Let $Q_n = P_2^n$, where $G^n$ denotes repeated application of the $\square$ operator, be the hypercube in $n$ dimensions.
Figure~\ref{fig:hypercubes} shows a DET:IC set for each of the hypercubes on $3 \le n \le 5$ dimensions.
From programmatic analysis, we believe these to be optimal DET:ICs.

\subsection{Cubic graphs}
\subsubsection{The infinite ladder}

\begin{figure}[ht]
    \centering
    \includegraphics[width=0.4\textwidth]{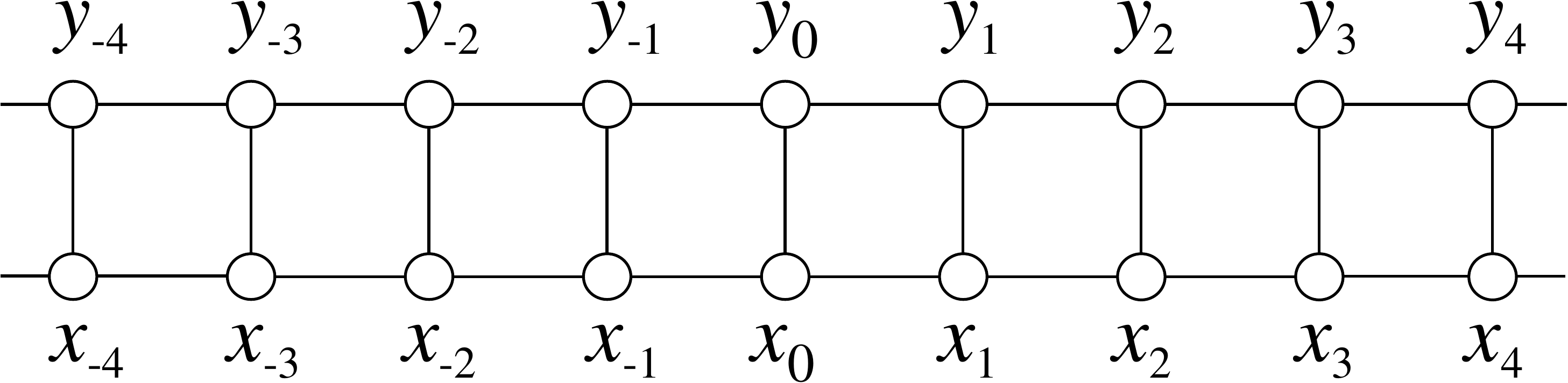}
    \caption{Ladder graph labeling scheme}
    \label{fig:ladder-labeling}
\end{figure}

\begin{figure}[ht]
    \centering
    \includegraphics[width=0.6\textwidth]{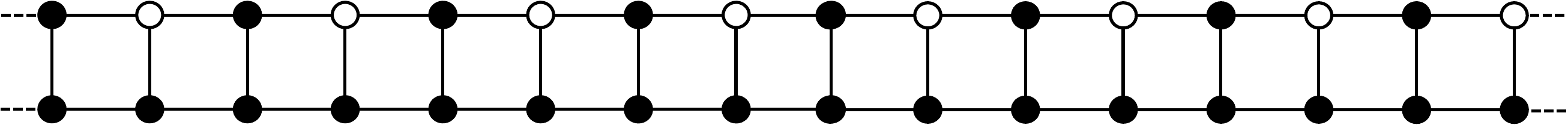}
    \caption{Unique optimal solution for the infinite ladder with $\textrm{DET:IC\%}(P_\infty \square P_2) = \frac{3}{4}$}
    \label{fig:ladder-det-ic-soln}
\end{figure}

\begin{theorem}
The infinite ladder graph has $\textrm{DET:IC\%}(P_\infty \square P_2) = \frac{3}{4}$.
\end{theorem}
\cbeginproof
The construction given by Figure~\ref{fig:ladder-det-ic-soln} is a density $\frac{3}{4}$ DET:IC on the infinite ladder graph.
We will prove that $\frac{3}{4}$ is the optimal value by showing that an arbitrary non-detector vertex can be associated with at least three detectors.
For $v \in V(G)$, let $R_5(v) = N(v) \cup \{u \in V(G) : |N(u) \cap N(v)| = 2\}$.
We impose that $x$ can be associated only with detector vertices within $R_5(x)$.
We will allow partial ownership of detectors, so a detector vertex, $v \in S$, contributes $\frac{1}{k}$, where $k = |R_5(v) \cap \overline{S}|$, toward the required total of three detectors.

Let $x_0 \notin S$ (see Figure~\ref{fig:ladder-labeling}).
If $y_0 \notin S$ then $x_1$ and $y_1$ cannot be distinguished, a contradiction.
Similarly, if $x_1 \notin S$ (or by symmetry $x_{-1} \notin S$), then $y_0$ and $y_1$ cannot be distinguished, a contradiction.
Therefore, by symmetry we see that no two non-detectors may be adjacent.
If $y_1 \notin S$ (or by symmetry $y_{-1} \notin S$) then $x_0$ and $y_1$ cannot be distinguished, a contradiction.
Thus, we see that $x_0 \notin S$ requires $\{x_{-1},x_1,y_{-1},y_0,y_1\} \subseteq S$, so there are five detectors in $R_5(x_0)$.
We will allow the possibility of a non-detector in $\{x_{-2},y_{-2}\}$ or $\{x_2,y_2\}$, but there may be only one in each pair since non-detectors cannot be adjacent.
Thus, $\{x_{-1},y_{-1}\}$ and $\{x_1,y_1\}$ could be associated with one non-detector other than $x_0$, but $y_0$ can only associate with $x_0$.
So $x_0$ has $\frac{2}{2} + \frac{2}{2} + \frac{1}{1} = 3$ associated detectors, completing the proof.
\cendproof

By Corollary~\ref{cor:detic-exist-alt-cubic} from Section 2, a cubic graph $G$ has a DET:IC if and only if $G$ is twin-free and triangle-free.
We introduce the notation $dom(v) = k$ to mean that $v$ is $k$-dominated.
We let $N^k(v) \subseteq V(G)$ denote the set of all vertices $x$ where there is a path of length $k$ from $v$ to $x$.
For each of the following propositions, we will assume that $S$ is a DET:IC for a cubic graph $G$.
\vspace{0.6em}

If $x \notin S$ and $xuvy$ is a path from $x$ to $y$, then $u$ and $v$ cannot be distinguished unless $y \in S$.
\begin{prop}\label{prop:detic-n3-s}
If $S \subseteq V(G)$ is a DET:IC for a cubic graph $G$ and $x \notin S$, then $N^3(x) \subseteq S$.
\end{prop}

Suppose $abcd$ is a 4-cycle in $G$. If $a,b \notin S$ we contradict Proposition~\ref{prop:detic-n3-s}. If $a,c \notin S$, then $a$ and $c$ are not distinguished. By symmetry, we have the following:
\begin{prop}\label{prop:detic-cubic-c4}
If $S$ is a DET:IC for a cubic graph $G$ and $abcd$ is a 4-cycle in $G$, then $|\{a,b,c,d\} \cap S| \ge 3$.
\end{prop}

If $xvy$ is a path from $x$ to $y$, then $v \in S$ is needed to 2-dominate $v$, and $(N(x) \cup N(y)) - \{v\} \subseteq N^3(x) \cup N^3(y)$, so Proposition~\ref{prop:detic-n3-s} yields the following result:
\begin{prop}\label{prop:detic-n2-s}
If $S \subseteq V(G)$ is a DET:IC for a cubic graph $G$ and $x,y \in V(G) - S$, with $y \in N^2(x)$, then $N(x) \cup N(y) \subseteq S$.
\end{prop}

If $xvy$ is a maximal path of detectors, then $x$ and $y$ are not distinguished, yielding the following:
\begin{prop}\label{prop:detic-cubic-max-p3}
If $S \subseteq V(G)$ is a DET:IC for cubic graph $G$ and $xvy$ is a path of detectors, then $dom(x) \ge 3$ or $dom(y) \ge 3$.
\end{prop}

\begin{definition}
Two vertices $p,q \in V(G)$ are called ``rivals" if there exists a 4-cycle $paqb$. And $p',q' \in (N(u) \cup N(v)) - \{p,a,q,b\}$ are called their ``friends".
\end{definition}

If $(x,y)$ are rivals with friends $(p,q)$, then $x,y \notin S$ is not allowed by Proposition~\ref{prop:detic-cubic-c4}, $x,q \notin S$ is not allowed by Proposition~\ref{prop:detic-n3-s}, and $p,q \notin S$ would cause $x$ and $y$ to not be distinguished. By symmetry, we have the following result:
\begin{prop}\label{prop:detic-cubic-rivals}
If $S$ is a DET:IC for cubic graph $G$, and $(x,y)$ are rivals with friends $(p,q)$, then $|\{x,y,p,q\} \cap S| \ge 3$.
\end{prop}

\begin{figure}[ht]
    \centering
    \begin{tabular}{c@{\hspace{2em}}||@{\hspace{2em}}c}
        \centered{\includegraphics[width=0.3\textwidth]{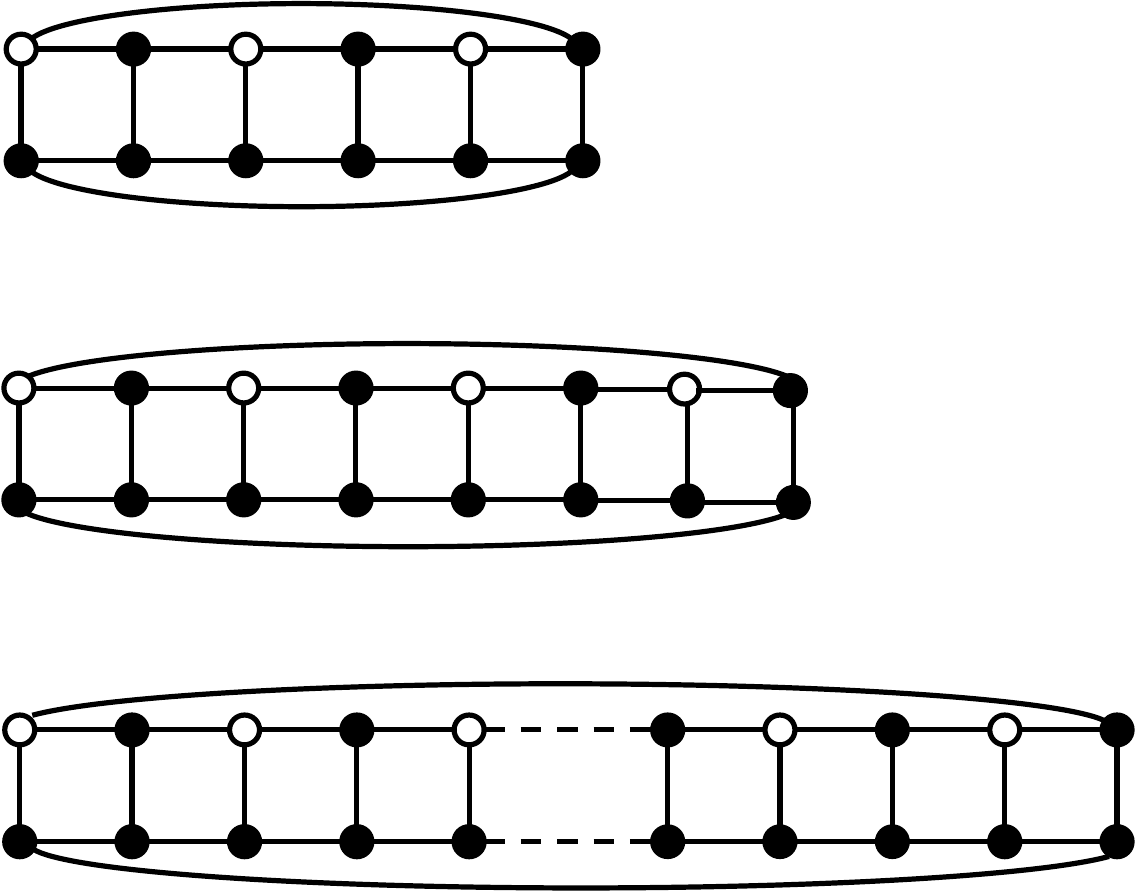}} &
        \centered{\includegraphics[width=0.5\textwidth]{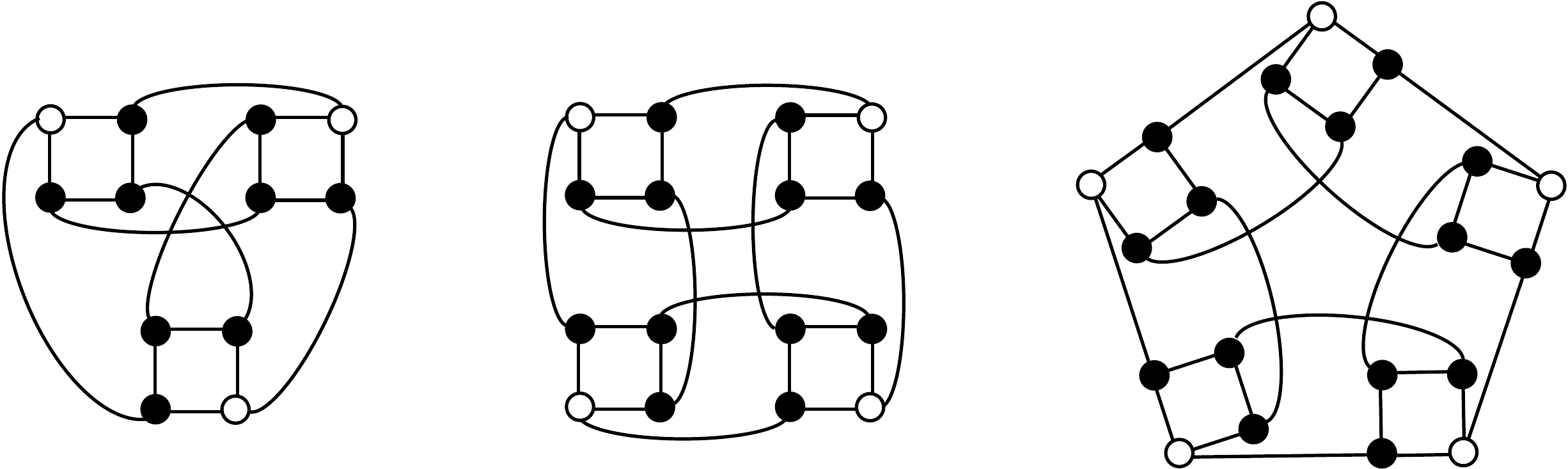}}
    \end{tabular}
    \caption{Two infinite families of cubic graphs with the low value of $\textrm{DET:IC\%}(G) = \frac{3}{4}$}
    \label{fig:detic-cubic-fams-3-4}
\end{figure}
\FloatBarrier

\subsubsection{Lower bound on DET:IC(G) for cubic}
The lowest value we have found for DET:IC in cubic graphs is $\ceil{\frac{3n}{4}}$.
Figures~\ref{fig:detic-cubic-fam-7-8} shows two infinite families of cubic graphs that achieve the density $\frac{3}{4}$.
We believe this density of $\frac{3}{4}$ is the minimum value for a cubic graph.

\begin{conjecture}
For a cubic graph $G$, $\textrm{DET:IC\%}(G) \ge \frac{3}{4}$.
\end{conjecture}

\subsubsection{Upper bound on DET:IC(G) for cubic}

\begin{figure}[ht]
    \centering
    \begin{tabular}{c@{\hskip 3em}c@{\hskip 3em}c@{\hskip 3em}c}
        \centered{\includegraphics[width=0.15\textwidth,angle=90]{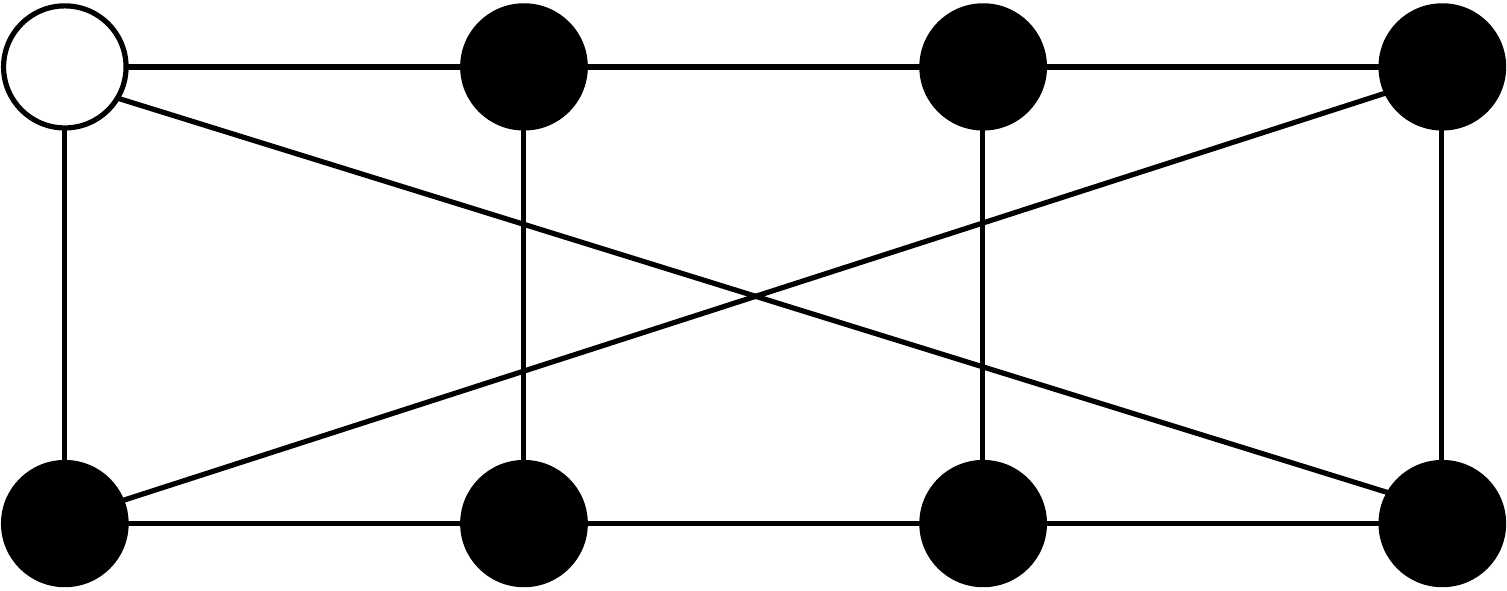}} &
        \centered{\includegraphics[width=0.18\textwidth,angle=90]{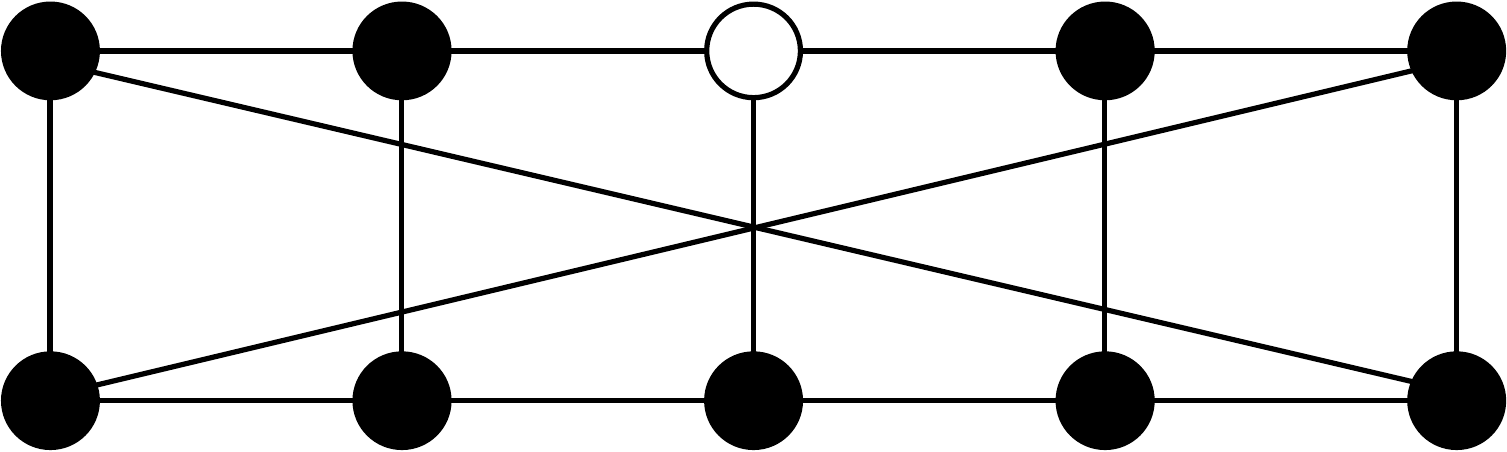}} &
        \centered{\includegraphics[width=0.2\textwidth]{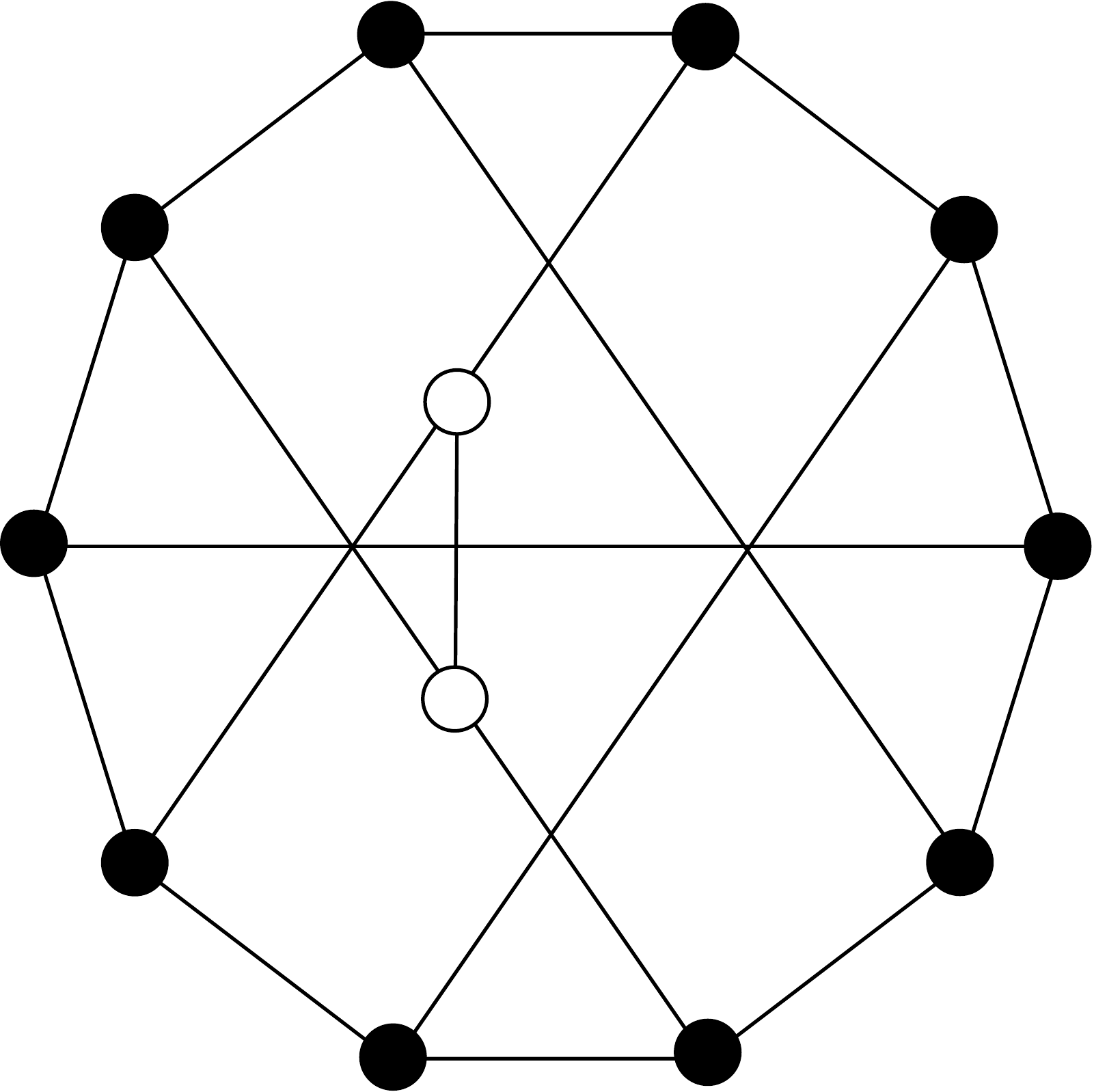}} &
        \centered{\includegraphics[width=0.2\textwidth]{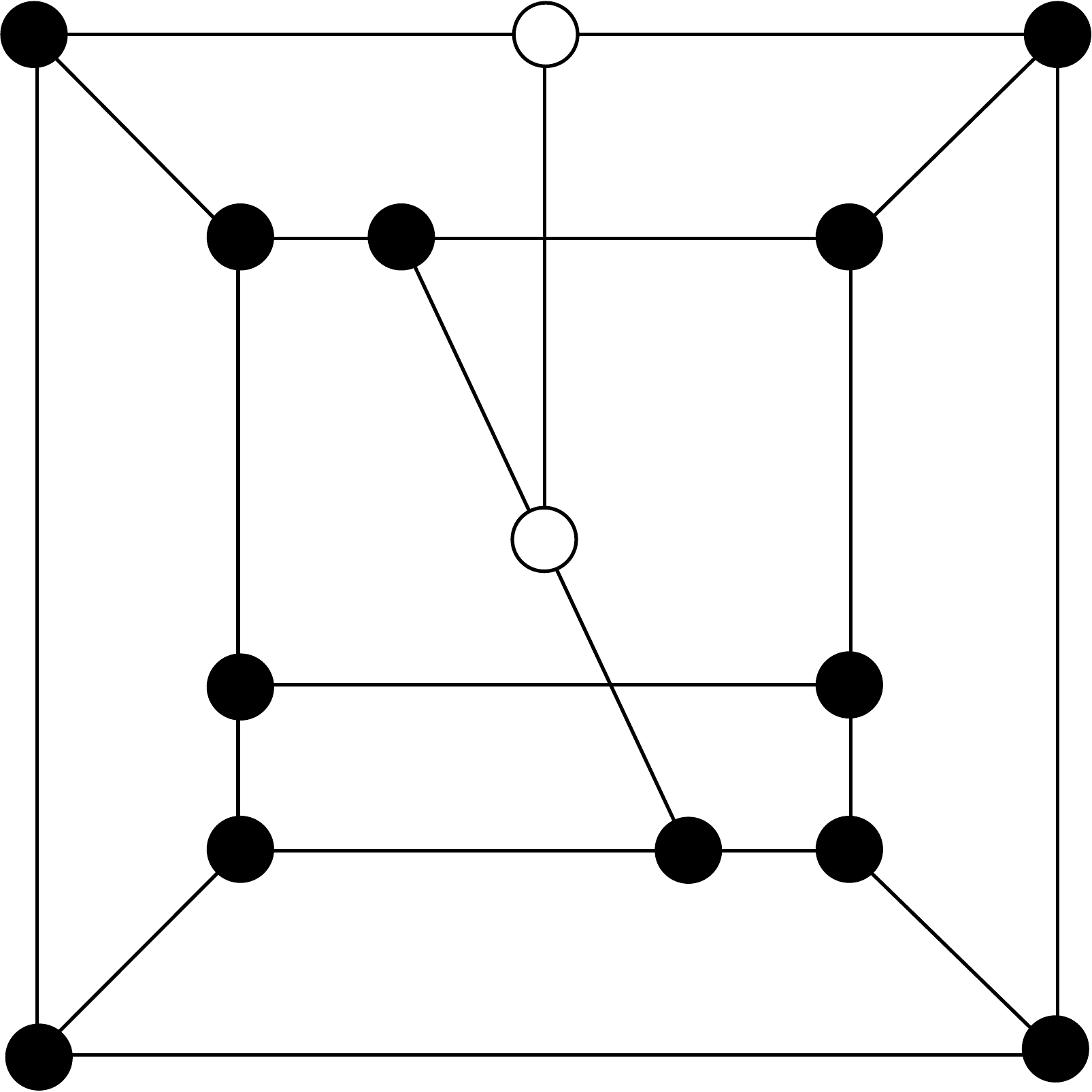}}\\
        $n=8$ & $n=10$ & $n=12$ & $n=14$ \\
        $\left(\frac{7}{8} = 0.875\right)$ & $\left(\frac{9}{10} = 0.9\right)$ & $\left(\frac{5}{6} \approx 0.8333\right)$ & $\left(\frac{6}{7} \approx 0.8571\right)$ \\\\
    \end{tabular}
    \begin{tabular}{c@{\hskip 2.5em}c@{\hskip 2.5em}c@{\hskip 2.5em}c}
        \centered{\includegraphics[width=0.19\textwidth]{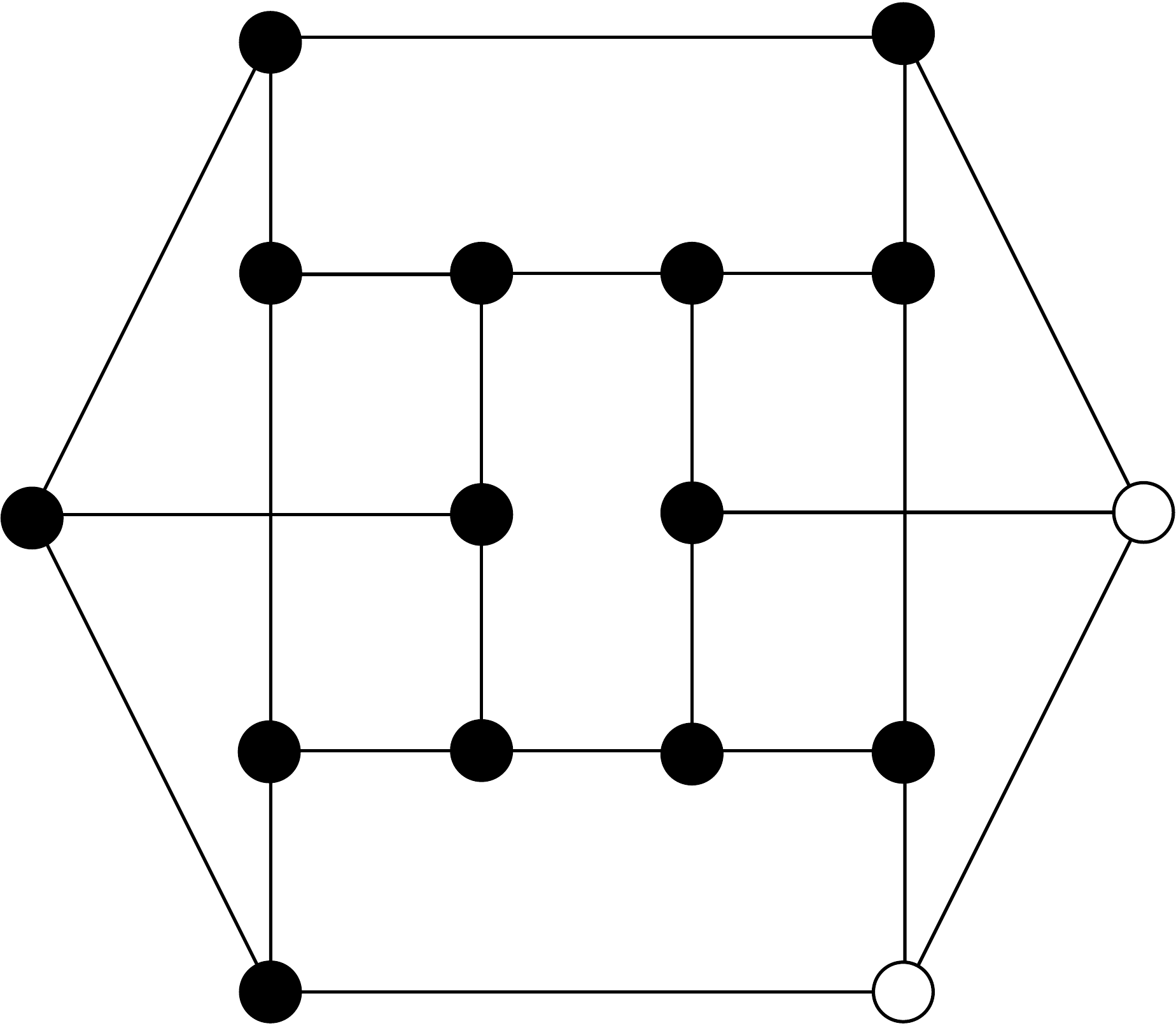}} &
        \centered{\includegraphics[width=0.23\textwidth,angle=90]{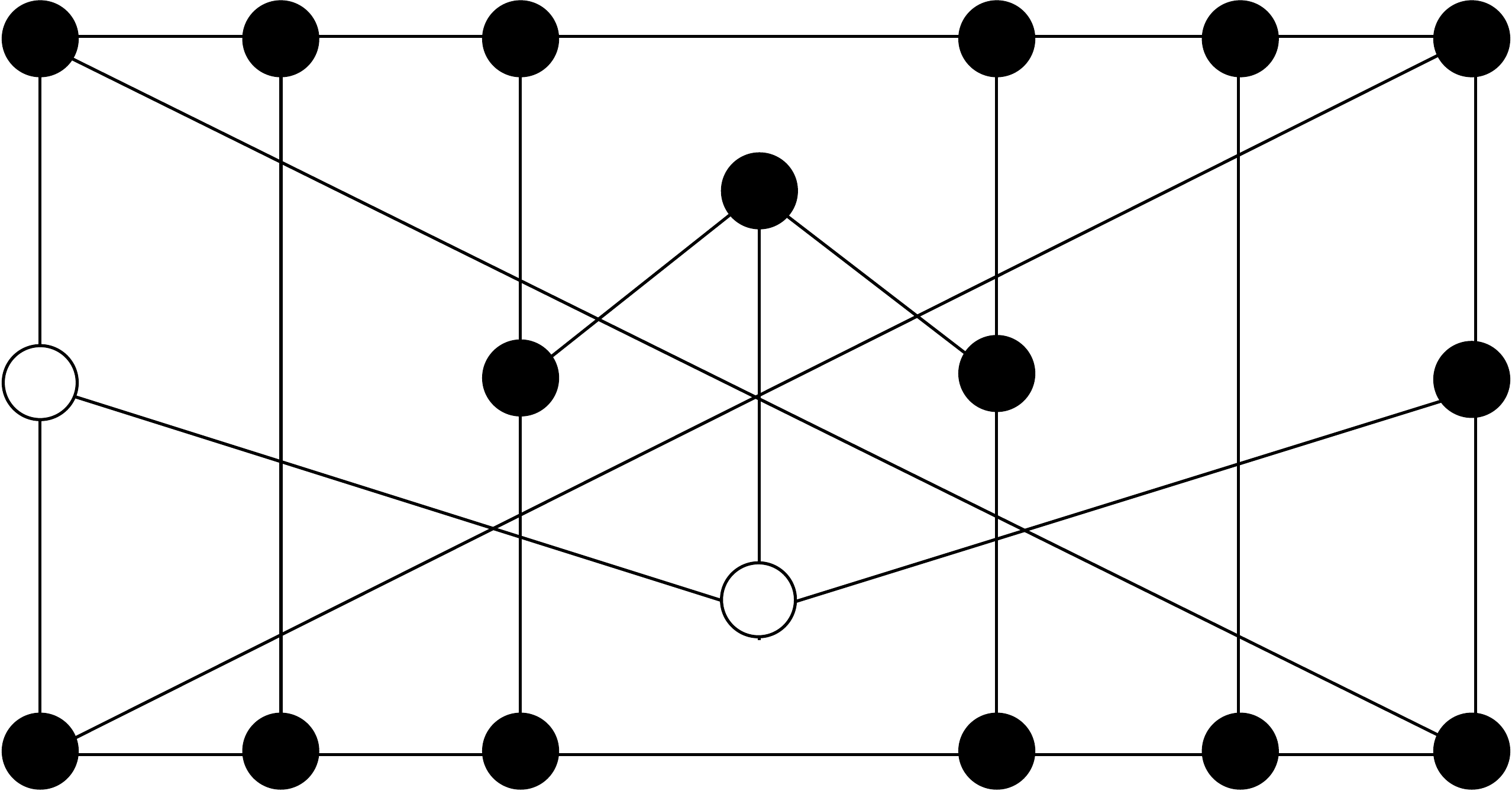}} &
        \centered{\includegraphics[width=0.28\textwidth,angle=90]{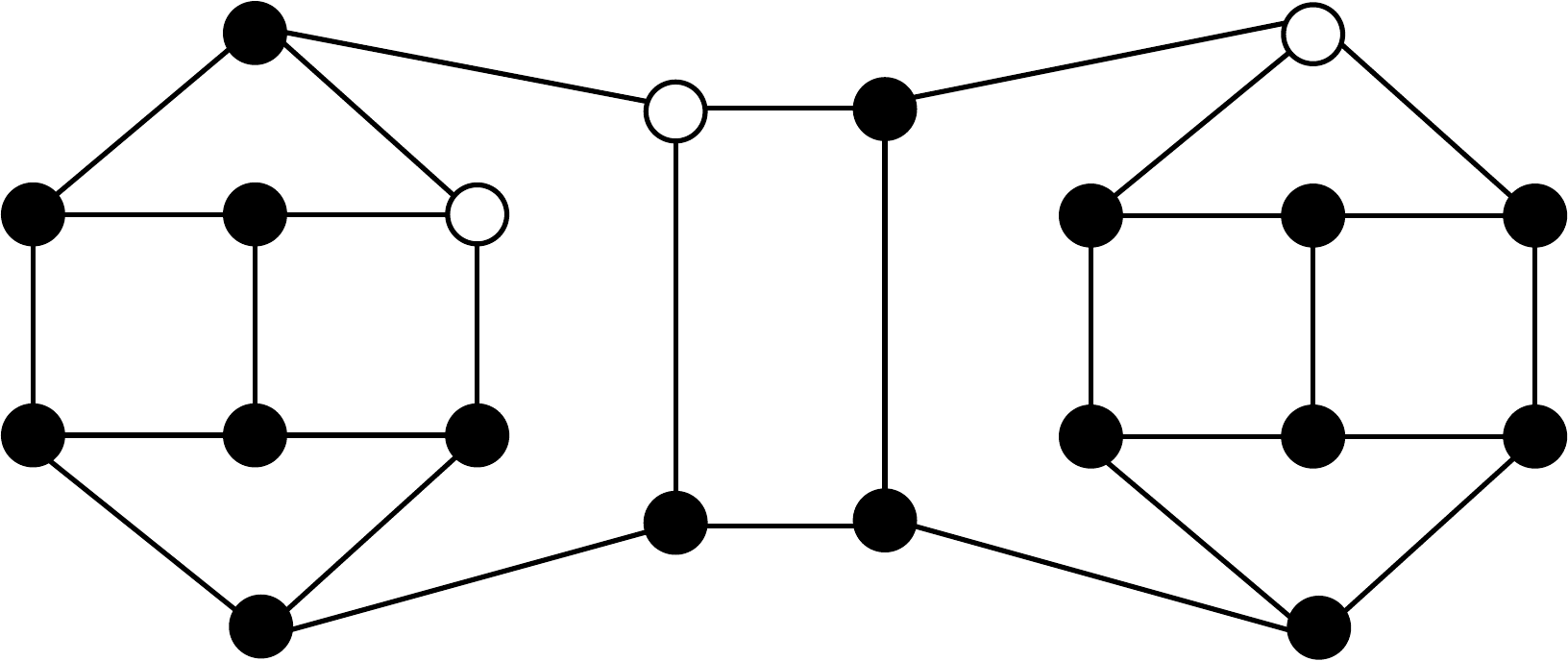}} &
        \centered{\includegraphics[width=0.22\textwidth]{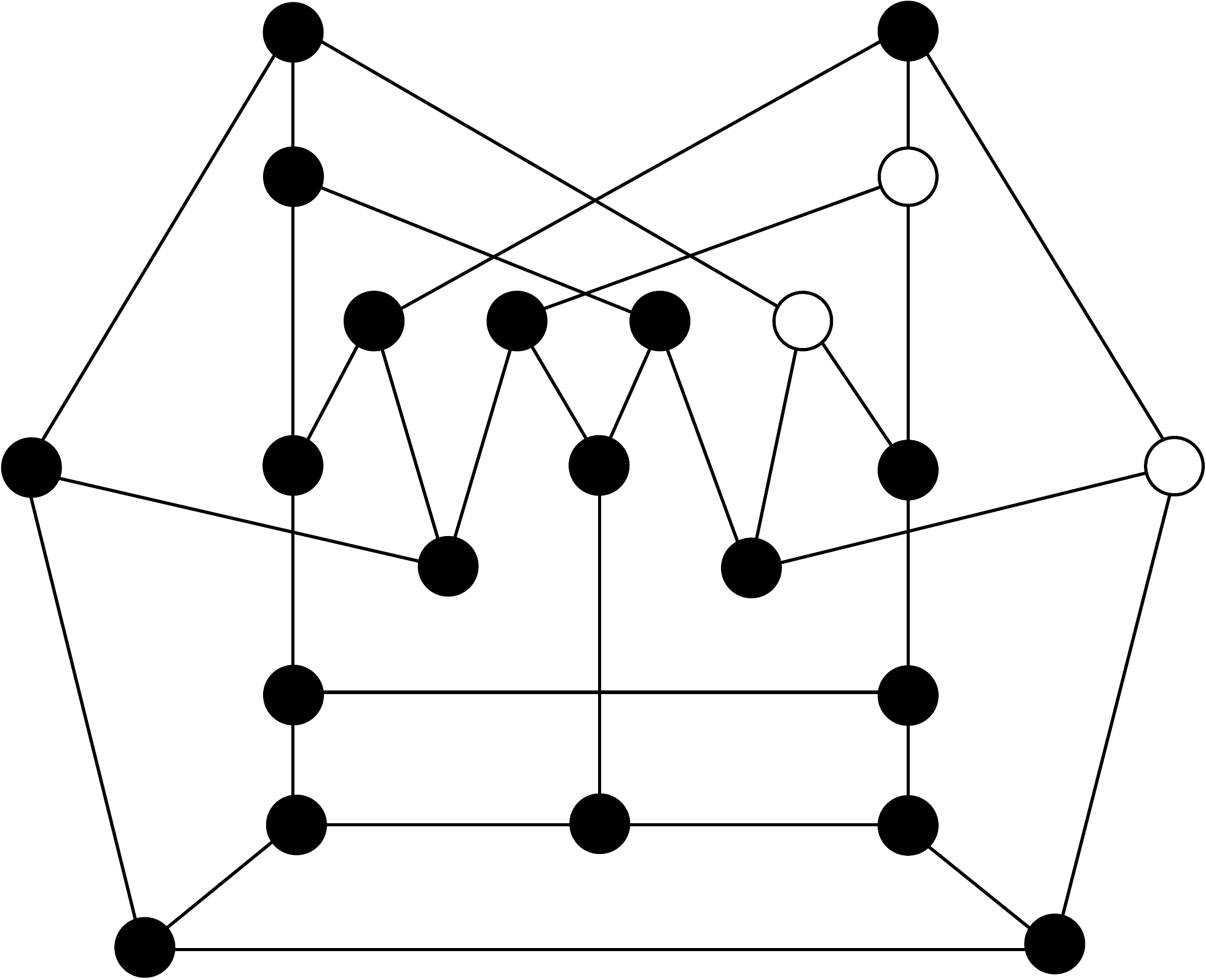}}\\
        $n=16$ & $n=18$ & $n=20$ & $n=22$ \\
        $\left(\frac{7}{8} = 0.875\right)$ & $\left(\frac{8}{9} \approx 0.8889\right)$ & $\left(\frac{17}{20} = 0.85\right)$ & $\left(\frac{19}{22} \approx 0.8636\right)$ \\\\    
    \end{tabular}
  \caption{Selected cubic graphs with highest $\textrm{DET:IC}(G)$ for $8 \le n \le 22$. Densities are shown in parenthesis.}
  \label{fig:cubic-fam-det-ic-ub}
\end{figure}

Figure~\ref{fig:cubic-fam-det-ic-ub} shows one example extremal cubic graph for each $n$ with $8 \le n \le 22$ with the highest value of $\textrm{DET:IC}(G)$.
The highest density we have found for DET:IC in cubic graphs is $\frac{9}{10}$, as shown in the $n=10$ graph from Figure~\ref{fig:cubic-fam-det-ic-ub}.
When $12 \le n \le 22$, the highest density is $\frac{8}{9}$, as shown in the $n=18$ graph; note that this is the only graph which achieves the density of $\frac{8}{9}$ when $n \le 22$.

\vspace{1.5em}
\FloatBarrier

\begin{figure}[ht]
    \centering
    \includegraphics[width=0.27\textwidth]{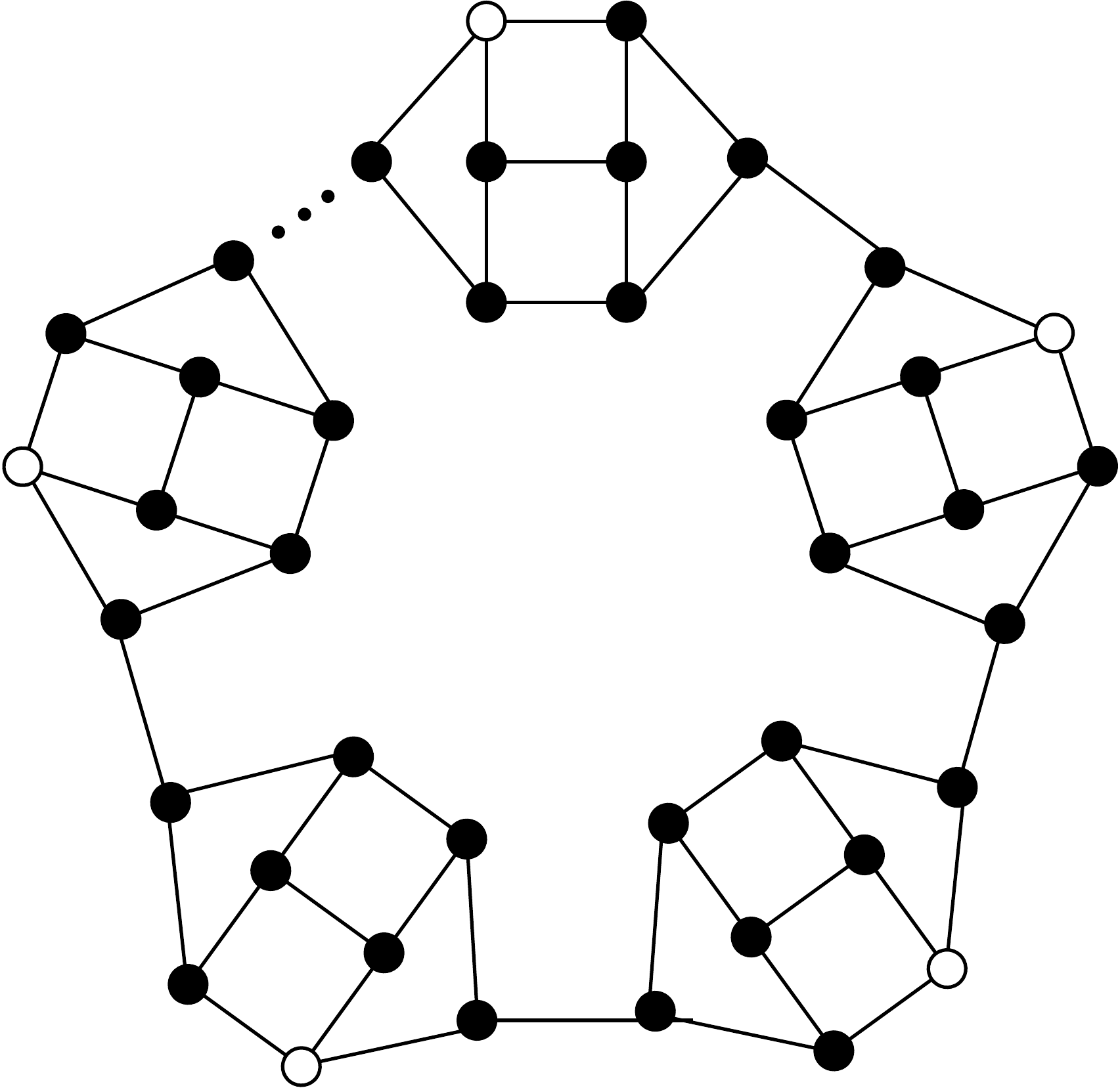}
    \caption{Infinite family of cubic graphs on $n=8k$ vertices with $\textrm{DET:IC\%}(G) = \frac{7}{8}$.}
    \label{fig:detic-cubic-fam-7-8}
\end{figure}

Figure~\ref{fig:detic-cubic-fam-7-8} shows an infinite family of cubic graphs that achieves the value $\frac{7}{8}$, which is proven in Theorem~\ref{theo:detic-cubic-fam-7-8}; this is the highest value among all infinite families of cubic graphs we have found thus far.

\begin{wrapfigure}[7]{r}{0.25\textwidth}
    \centering
    \includegraphics[width=0.22\textwidth]{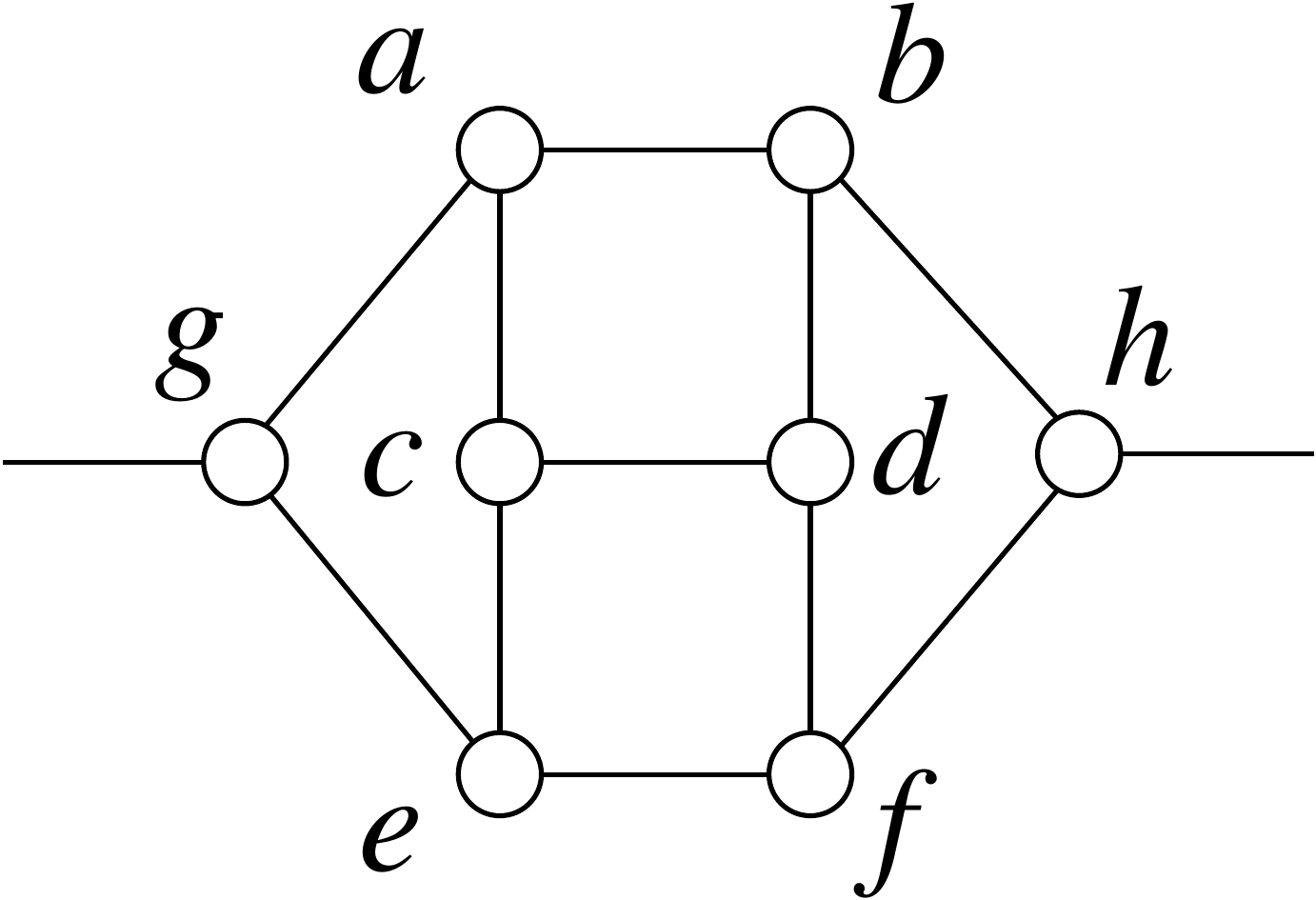}
    \caption{Subgraph $H$}
    \label{fig:detic-cubic-fam-subgraph}
\end{wrapfigure}

\begin{lemma}\label{lem:detic-cubic-fam-7-8}
If $S$ is a DET:IC for a cubic graph containing subgraph $H$, as shown in Figure~\ref{fig:detic-cubic-fam-subgraph}, then $|V(H) \cap S| \ge 7$.
\end{lemma}
\cbeginproof
Let $H$ be the subgraph on 8 vertices shown in Figure~\ref{fig:detic-cubic-fam-subgraph}; note that the two outgoing edges from $g$ and $h$ can go anywhere in the graph other than the vertices of $H$, which are already degree 3.

If $c \notin S$, then Proposition~\ref{prop:detic-cubic-c4} yields that $\{a,b,d,e,f,g\} \subseteq S$, and Proposition~\ref{prop:detic-n3-s} yields that $h \in S$ so there is at most one non-detector in $H$.
Otherwise, we assume that $c \in S$, and $d \in S$ by symmetry.
If $g \notin S$, then Proposition~\ref{prop:detic-cubic-c4} yields that $\{a,e\} \subseteq S$, Proposition~\ref{prop:detic-n3-s} yields that $h \in S$, and Proposition~\ref{prop:detic-cubic-rivals} applied to rival pairs $(a,d)$ and $(d,e)$ yields that $f \in S$ and $b \in S$, respectively, and we are done.
Otherwise, we assume that $g \in S$ and by symmetry $h \in S$.
We now see that applying Proposition~\ref{prop:detic-cubic-rivals} to rival pair $(a,e)$ yields that there may be only one non-detector in the remaining unknown vertices, $\{a,b,e,f\}$, completing the proof.
\cendproof

\FloatBarrier

\begin{theorem}\label{theo:detic-cubic-fam-7-8}
The infinite family of cubic graphs given in Figure~\ref{fig:detic-cubic-fam-7-8} has value $\textrm{DET:IC\%}(G) = \frac{7}{8}$.
\end{theorem}
\begin{proof}
Each $G$ in the family on $n=8k$ vertices is composed (exclusively) of $k \ge 1$ copies of Subgraph $H$ from Figure~\ref{fig:detic-cubic-fam-subgraph}.
Thus, Lemma~\ref{lem:detic-cubic-fam-7-8} gives us that there are at most $k$ non-detectors.
It can be shown that $S = V(G) - A$ is a DET:IC for $G$, where $A = \{a_1,\hdots,a_k\}$ is the set of $a$ vertices from the $k$ copies of subgraph $H$; this solution is shown in Figure~\ref{fig:detic-cubic-fam-7-8}.
Thus, we achieve the value of $\frac{7}{8}$, completing the proof.
\end{proof}

\subsection{Infinite Grids}
\begin{figure}[ht]
    \centering
    \begin{tabular}{c@{\hskip 4em}c}
        \includegraphics[width=0.4\textwidth]{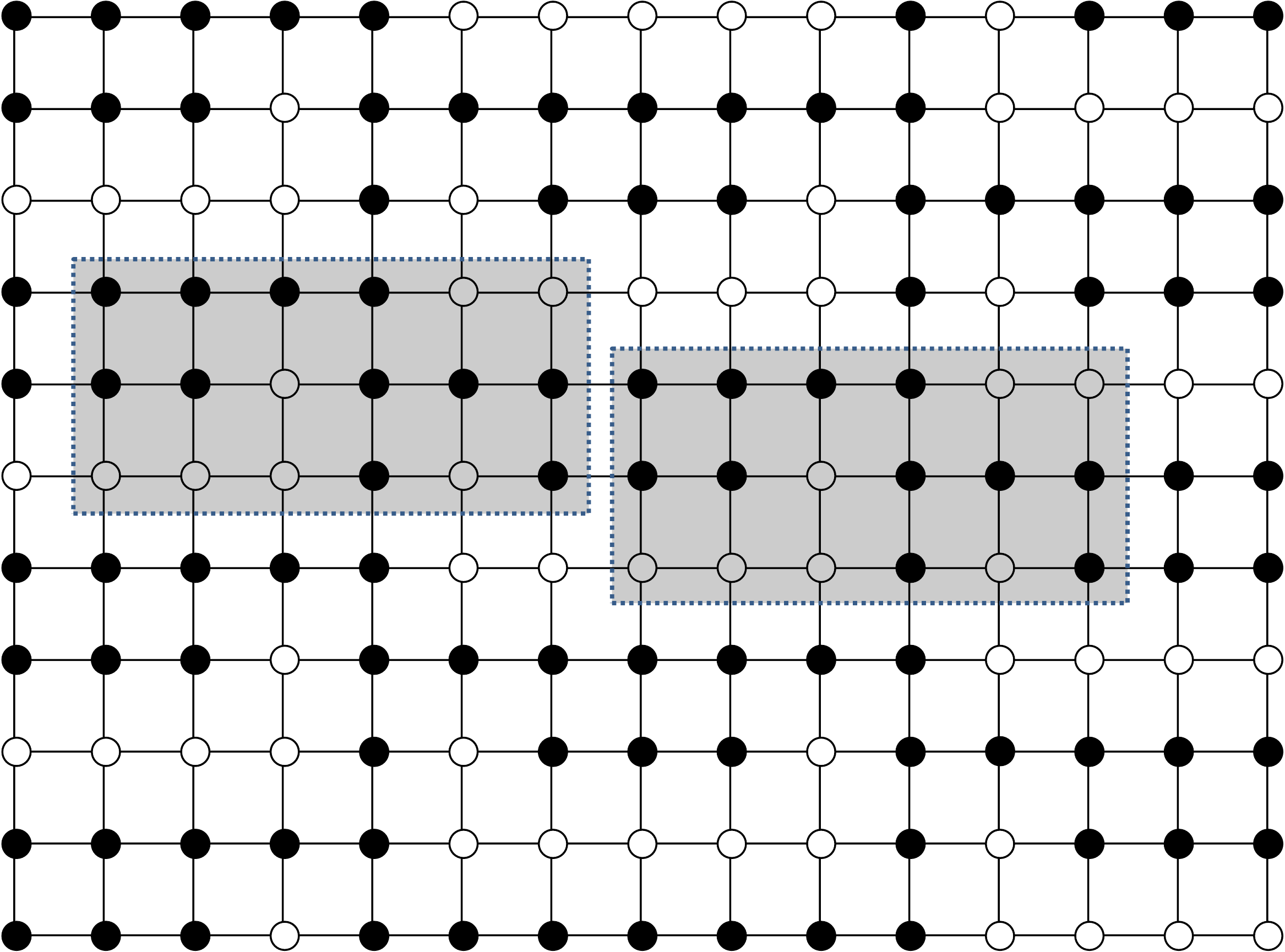} & \includegraphics[width=0.4\textwidth]{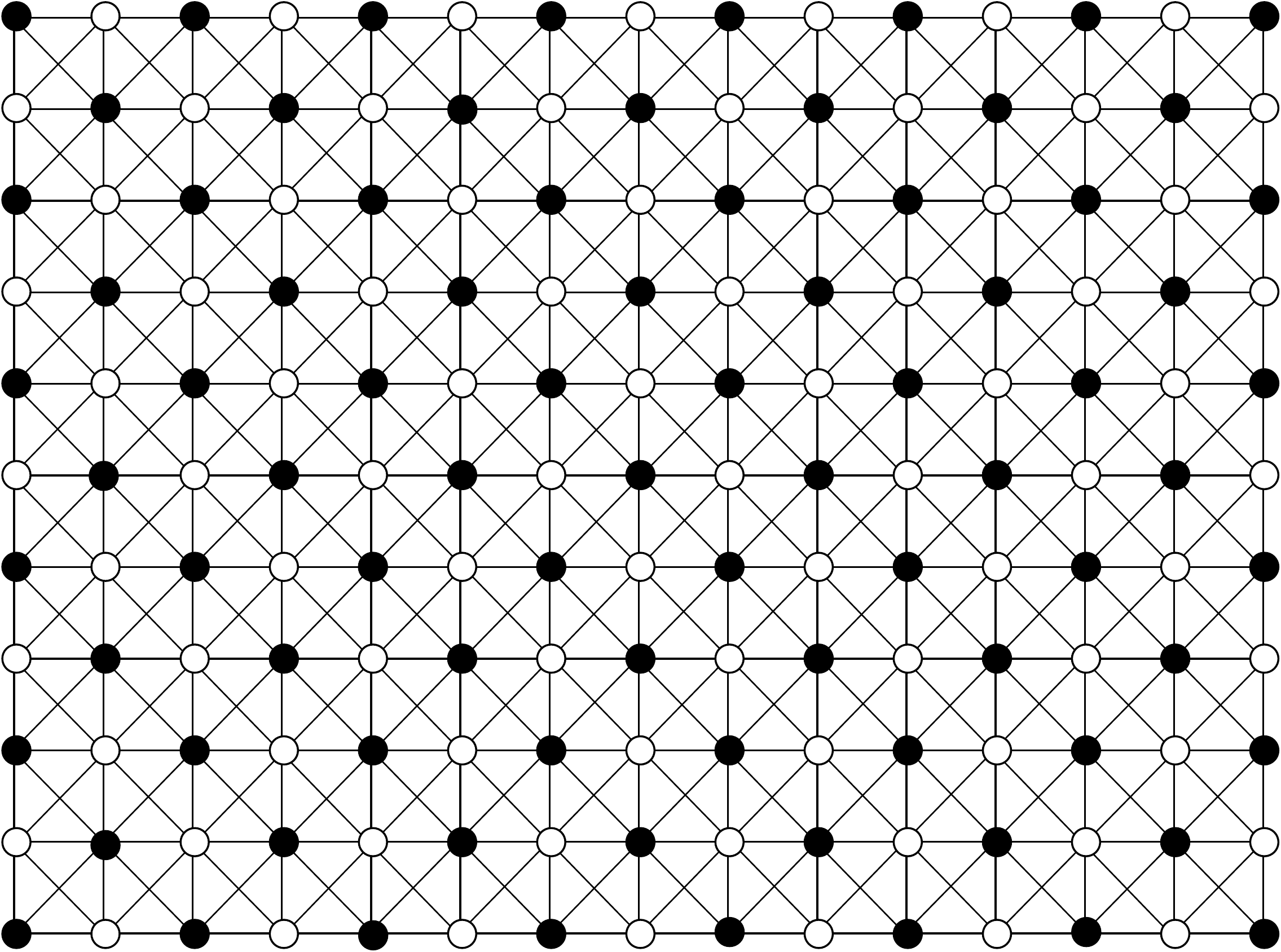} \\ (a) & (b) \\ \\
        \includegraphics[width=0.4\textwidth]{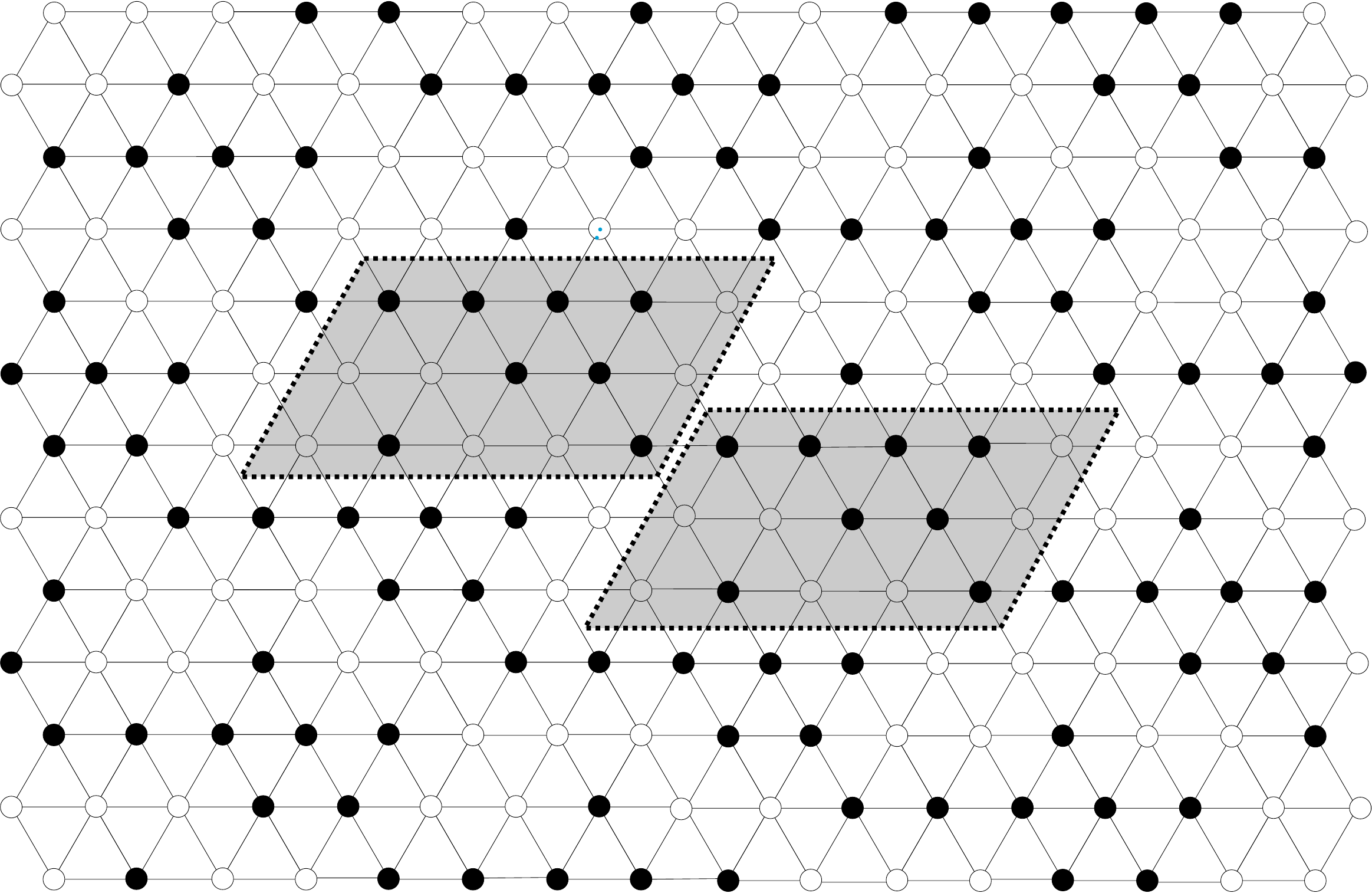} & \includegraphics[width=0.4\textwidth]{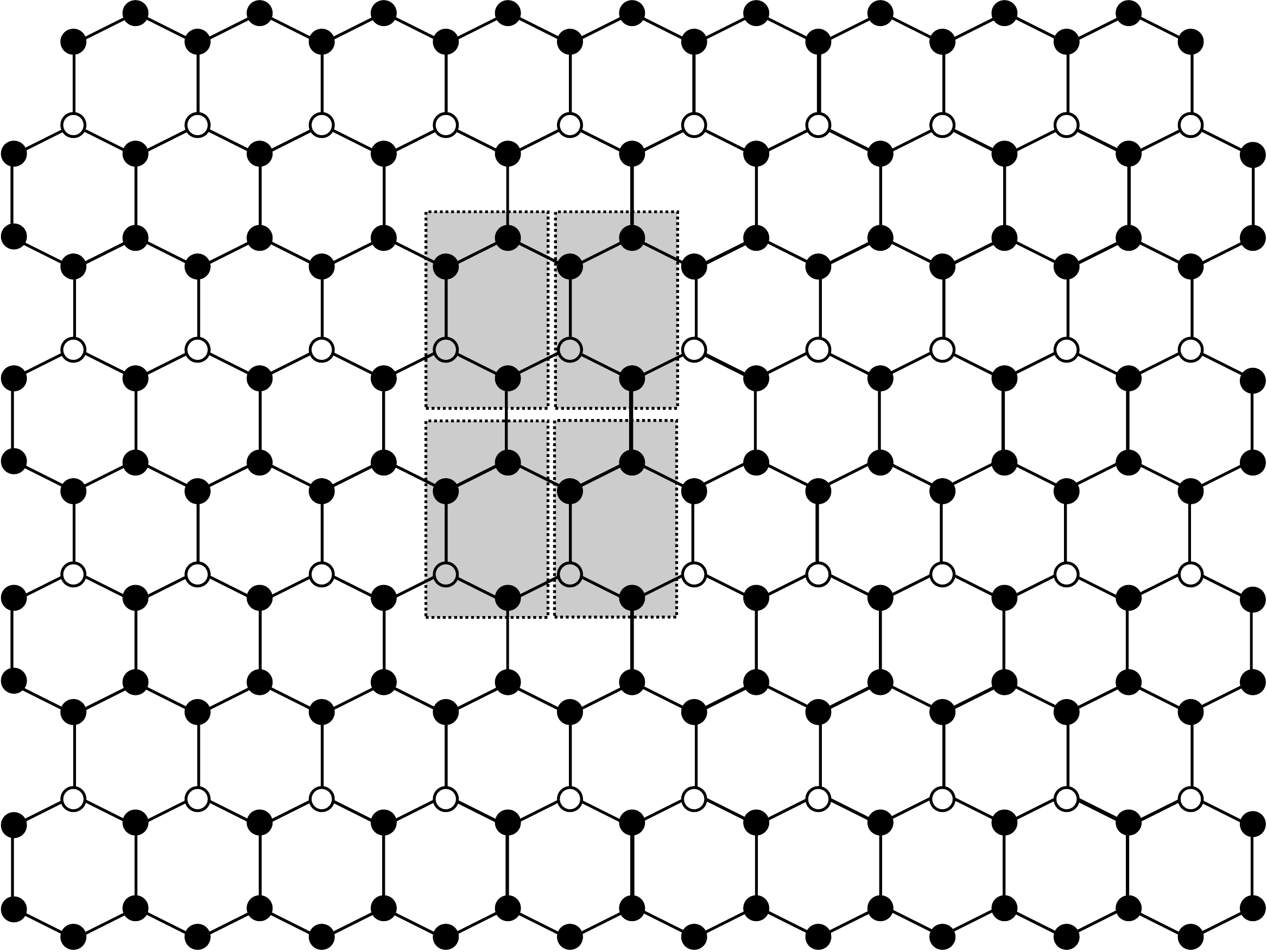} \\ (c) & (d)
    \end{tabular}
    \caption{Our best constructions of DET:IC on SQR (a), KNG (b), TRI (c), and HEX (d). Shaded vertices denote detectors.}
    \label{fig:inf-grids-det-ic-solns}
\end{figure}

We have the following theorem on the results on DET:IC for some infinite grids, and the solutions that achieve each of the upper bounds can be found in Figure~\ref{fig:inf-grids-det-ic-solns}.

\begin{theorem}
The upper and lower bounds on DET:IC:
\begin{enumerate}[label=\roman*]
    \item For the infinite hexagonal grid HEX, $\frac{12}{17} \le \textrm{DET:IC\%}(HEX) \le \frac{3}{4}$.
    \item For the infinite square grid SQR, $\frac{10}{17} \le \textrm{DET:IC\%}(SQR) \le \frac{11}{18}$.
    \item For the infinite triangular grid TRI, $\frac{30}{61} \le \textrm{DET:IC\%}(TRI) \le \frac{8}{15}$.
    \item For the infinite king grid KNG,
$\frac{40}{99} \le \textrm{DET:IC\%}(KNG) \le \frac{1}{2}$.
\end{enumerate}
\end{theorem}

\FloatBarrier

\bibliographystyle{acm}
\bibliography{refs}

\end{document}